\documentclass[a4paper,10pt]{article}
\usepackage[utf8x]{inputenc}
\usepackage{tracefnt,amsmath,tabu,array}
\usepackage{amssymb,graphicx,setspace,amsfonts,amsbsy}
\usepackage{pifont,latexsym,ifthen,amsthm,rotating,calc,textcase,booktabs,cancel,slashed}

\newtheorem{theorem}{Theorem}[section]
\newtheorem{lemma}[theorem]{Lemma}
\newtheorem{corollary}[theorem]{Corollary}
\newtheorem{proposition}[theorem]{Proposition}

\newtheorem{remark}[theorem]{Remark}
\newcommand{\filledbox}{\leavevmode
  \hbox to.77778em{%
  \hfil\vbox to.675em{\hrule width.6em height.6em}\hfil}}

\newcommand{\Rm}{{\mathbb R}}

\newcommand{\eps}{\varepsilon}

\begin{document}
\tabulinesep=1.0mm
\title{Explicit formula of radiation fields of free waves with applications on channel of energy}

\author{Liang Li, Ruipeng Shen and Lijuan Wei\\
Centre for Applied Mathematics\\
Tianjin University\\
Tianjin, China
}

\maketitle

\begin{abstract}
 In this work we give a few explicit formulas regarding the radiation fields of linear free waves. We then apply these formulas on the channel of energy theory. We characterize all the radial weakly non-radiative solutions in all dimensions and give a few new exterior energy estimates. 
\end{abstract}

\section{Introduction}

\subsection{Background and topics} 
The semi-linear wave equation
\[
 \partial_t^2 u - \Delta u = \pm |u|^{p-1} u,\qquad (x,t) \in \Rm^d \times \Rm;
\]
especially the energy critical case $p = 1 + 4/(d-2)$, has been extensively studied by many mathematicians in the past few decades. Please see, for example, Kapitanski \cite{loc1} and Lindblad-Sogge \cite{ls} for local existence and well-posedness; Ginibre-Soffer-Velo \cite{locad1}, Grillakis \cite{mg1,mg2}, Kenig-Merle \cite{kenig}, Nakanishi \cite{enscatter1, enscatter2} and Shatah-Struwe \cite{ss1, ss2} for global existence, regularity, scattering and blow-up. Since the semi-linear wave equation can be viewed as a small perturbation of the homogenous linear wave equation in many situations, especially when we consider the asymptotic behaviours of solutions as spatial variables or time tends to infinity, it is important to first understand the asymptotic behaviours of solutions to the homogenous linear wave equation, i.e. free waves. This work is concerned with two important tools to understand the asymptotic behaviours of free waves: radiation field and channel of energy. We first introduce a few necessary notations. Throughout this work we consider the homogenous linear wave equation with initial data in the energy space
\begin{equation} \label{LW}
 \left\{\begin{array}{ll} \partial_t^2 u - \Delta u = 0, & (x,t) \in \Rm^d \times \Rm; \\
 u|_{t=0} = u_0 \in \dot{H}^1(\Rm^d); \\
 u_t|_{t=0} = u_1 \in L^2(\Rm^d). \end{array}\right.
\end{equation}
In this work we also use the notation $\mathbf{S}_L(u_0,u_1)$ to represent the free wave $u$ defined above. If it is necessary to mention the velocity $u_t$, we use the notation 
\[
 \mathbf{S}_L(t) \begin{pmatrix} u_0\\ u_1 \end{pmatrix} = \begin{pmatrix} u(\cdot,t)\\ u_t(\cdot, t)\end{pmatrix} \in \dot{H}^1 \times L^2. 
\]
It is well known that the linear wave propagation preserves the $\dot{H}^1 \times L^2$ norm, i.e. the energy conservation law holds. ($\nabla_{x,t} u = (\nabla u, u_t)$)
\[
 \int_{\Rm^d} |\nabla_{x,t} u(x,t)|^2 dx = \int_{\Rm^d} (|\nabla u_0|^2 + |u_1|^2) dx. 
\]
Now we make a brief review of radiation field and channel of energy method. 

\paragraph{Radiation field} The asymptotic behaviour of free waves at the energy level can be characterized by the following theorem.

\begin{theorem}[Radiation filed] \label{radiation}
Assume that $d\geq 3$ and let $u$ be a solution to the free wave equation $\partial_t^2 u - \Delta u = 0$ with initial data $(u_0,u_1) \in \dot{H}^1 \times L^2(\Rm^d)$. Then 
\[
 \lim_{t\rightarrow \pm \infty} \int_{\Rm^d} \left(|\nabla u(x,t)|^2 - |u_r(x,t)|^2 + \frac{|u(x,t)|^2}{|x|^2}\right) dx = 0
\]
 and there exist two functions $G_\pm \in L^2(\Rm \times \mathbb{S}^{d-1})$ so that
\begin{align*}
 \lim_{t\rightarrow \pm\infty} \int_0^\infty \int_{\mathbb{S}^{d-1}} \left|r^{\frac{d-1}{2}} \partial_t u(r\theta, t) - G_\pm (r+t, \theta)\right|^2 d\theta dr &= 0;\\
 \lim_{t\rightarrow \pm\infty} \int_0^\infty \int_{\mathbb{S}^{d-1}} \left|r^{\frac{d-1}{2}} \partial_r u(r\theta, t) \pm G_\pm (r+t, \theta)\right|^2 d\theta dr & = 0.
\end{align*}
In addition, the maps $(u_0,u_1) \rightarrow \sqrt{2} G_\pm$ are a bijective isometries form $\dot{H}^1 \times L^2(\Rm^d)$ to $L^2 (\Rm \times \mathbb{S}^{d-1})$.
\end{theorem}

\noindent This has been known for more than 50 years, at least in the 3-dimensional case. Please see Friedlander \cite{radiation1, radiation2}, for example. The version of radiation field theorem given above and a proof for all dimensions $d\geq 3$ can be found in Duyckaerts et al. \cite{dkm3}. In addition, there is also a 2-dimensional version of radiation field theorem. The statement in dimension $2$ can be given in almost the same way as in the higher dimensional case, except that the limit 
\[ 
 \lim_{t\rightarrow \pm \infty} \int_{\Rm^2} \frac{|u(x,t)|^2}{|x|^2} dx = 0
\]
no longer holds. A proof by Radon transform for all dimensions $d\geq 2$ can be found in Katayama \cite{katayamaradiation}, where the statement of the theorem is slightly different. Throughout this work we call the function $G_\pm$ radiation profiles and use the notations $\mathbf{T}_\pm$ for the linear map $(u_0,u_1)\rightarrow G_\pm$.

\paragraph{Channel of energy} The second tool is the channel of energy method, which plays an important role in the study of wave equation in the past decade. This method is first introduced in 3-dimensional case by Duyckaerts-Kenig-Merle \cite{tkm1} and then in 5-dimensional case by Kenig-Lawrie-Schlag \cite{channel5d}. This method was used in the proof of solition resolution conjecture of energy critical wave equation with radial data in all odd dimensions $d\geq 3$ by Duyckaerts-Kenig-Merle \cite{se, oddhigh}. It can also be used to show the non-existence of minimal blow-up solutions in a compactness-rigidity argument in the energy super or sub-critical case. Please see, for example, Duyckaerts-Kenig-Merle \cite{dkm2} and Shen \cite{shen2}.  Roughly speaking, the channel of energy method discusses the amount of energy located in an exterior region as time tends to infinity:
\[
 \lim_{t\rightarrow \pm \infty} \int_{|x|>|t|+R} |\nabla_{x,t} u(x,t)|^2 dx.
\]
Here the constant $R\geq 0$. Since the energy travels at a finite speed, the energy in the exterior region $\{x: |x|>|t|+R\}$ decays as $|t|$ increases. Thus the limits above are always well-defined. We can also give the exact value of the limit in term of the radiation field:
\begin{equation} \label{relationship radiation and channel}
 \lim_{t\rightarrow \pm \infty} \int_{|x|>|t|+R} |\nabla_{x,t} u(x,t)|^2 dx = 2\int_R^\infty \int_{\mathbb{S}^{d-1}} |G_\pm (s, \theta)|^2 d\theta ds. 
\end{equation}
We first introduce a few already known results. We start with the odd dimensions. 
\begin{proposition} [see Duyckaerts-Kenig-Merle \cite{dkmnonradial}] \label{situation R0}
Assume that $d\geq 3$ is an odd integer. All solutions to $\partial_t^2 u - \Delta u = 0$ satisfies 
\begin{equation} \label{radiative identity}
  \sum_{\pm} \lim_{t\rightarrow \pm \infty} \int_{|x|>|t|} |\nabla_{x,t} u(x,t)|^2 dx= \int_{\Rm^d} |\nabla_{x,t} u(x,0)|^2 dx. 
\end{equation}
\end{proposition}
\noindent As a result, we have 
\begin{corollary} \label{non radiative is zero}
Assume that $d\geq 3$ is odd. Then $u\equiv 0$ is the only free wave $u$ satisfying 
\[
 \lim_{t\rightarrow \pm \infty} \int_{|x|>|t|} |\nabla_{x,t} u(x,t)|^2 dx = 0.
\]
\end{corollary}
\noindent In the contrast, if $R>0$, the subspace of $\dot{H}^1 \times L^2(\Rm^d)$ defined by
\begin{equation} \label{def of PR}
 P(R) = \left\{(u_0,u_1)\in \dot{H}^1\times L^2(\Rm^d): \lim_{t\rightarrow \pm \infty} \int_{|x|>R+|t|} \left|\nabla_{t,x} \mathbf{S}_L(u_0,u_1)(x,t)\right|^2 dx = 0\right\}. 
\end{equation} 
does contain initial data $(u_0,u_1)$ whose support is essentially bigger than $\{x: |x|\leq R\}$. The free waves $u$ satisfying
\[
 \lim_{t\rightarrow \pm \infty} \int_{|x|>R+|t|} \left|\nabla_{t,x} u(x,t)\right|^2 dx = 0
\]
 are usually called ($R$-weakly) non-radiative solutions. If the dimension is odd, these solutions are well-understood in the radical case: 
\begin{theorem} [See Kenig et al \cite{channel}, the proof uses radial Fourier transform]\label{radialpi} 
In any odd dimension $d\geq 1$, every radial solution $u$ to \eqref{LW} satisfies 
\begin{equation} \label{radial exterior energy}
 \max_{\pm} \lim_{t\rightarrow \pm \infty} \int_{r>|t|+R} |\nabla_{x,t} u(r,t)|^2 r^{d-1} dr\geq \frac{1}{2} \left\|\Pi_{P_{rad}(R)}^\perp (u_0,u_1)\right\|_{\dot{H}^1\times L^2(r\geq R: r^{d-1}dr)}.
\end{equation}
Here 
\[
 P_{rad}(R) \doteq \hbox{Span} \left\{(r^{2k_1-d},0), (0, r^{2k_2-d}): k_1, k_2\in \mathbb{N}; 1\leq k_1 \leq \frac{d+2}{4}, 1\leq k_2 \leq \frac{d}{4}\right\}.
\]
$\Pi_{P_{rad}(R)}^\perp$ is the orthogonal projection from $\dot{H}^1 \times L^2(r\geq R: r^{d-1} dr)$ onto the complement of the finite-dimensional subspace $P_{rad}(R)$. 
\end{theorem}
 
\noindent The case of even dimensions is much more complicated and subtle. C\^{o}te-Kenig-Schlag \cite{channeleven} shows that in general the inequality
\[
  \sum_{\pm} \lim_{t\rightarrow \pm \infty} \int_{|x|>|t|} |\nabla_{x,t} u(x,t)|^2 dx \geq C \int_{\Rm^d} |\nabla_{x,t} u(x,0)|^2 dx
\]
does not hold for any positive constant $C$ in even dimensions. But a similar inequality holds in the radial case for either initial data $(u_0,0)$, if $d = 0\, \hbox{mod}\, 4$, or $(0,u_1)$, if $d= 2\, \hbox{mod}\, 4$. More precisely we have
\begin{align}
 \lim_{t\rightarrow \pm \infty} \int_{|x|>|t|} |\nabla_{x,t} \mathbf{S}_L(u_0,0)(x,t)|^2 dx & \geq \frac{1}{2} \int_{\Rm^d} |\nabla u_0(x)|^2 dx, \quad d=4k, \, k\in \mathbb{N}; \label{exterior full 4k} \\
  \lim_{t\rightarrow \pm \infty} \int_{|x|>|t|} |\nabla_{x,t} \mathbf{S}_L(0,u_1)(x,t)|^2 dx & \geq \frac{1}{2} \int_{\Rm^d} |u_1(x)|^2 dx, \quad d=4k+2, \, k\in \mathbb{N} \label{exterior full 4k plus 2}.
\end{align}
In addition, Duychaerts-Kenig-Merle \cite{oddtool} shows that the only non-radiative solution is still zero solution in even dimensions $d\geq 4$, i.e. Corollary \ref{non radiative is zero} still holds for even dimensions $d\geq 4$, even in the non-radial case. Much less is known about the exterior energy estimate in the region $\{x: |x|>R+|t|\}$ with $R>0$. Dyuchaerts at el. \cite{coratational} proves the exterior energy estimate in dimension 4 and 6 if the initial data are radial:
\begin{align*}
 \lim_{t\rightarrow \pm \infty} \int_{|x|>|t|+R} |\nabla_{x,t} \mathbf{S}_L(u_0,0)(x,t)|^2 dx & \geq \frac{3}{10}  \|\Pi_R^\perp u_0\|_{\dot{H}^1(\{x\in \Rm^4: |x|>R\})}^2;\\
  \lim_{t\rightarrow \pm \infty} \int_{|x|>|t|+R} |\nabla_{x,t} \mathbf{S}_L(0,u_1)(x,t)|^2 dx & \geq \frac{3}{10}  \|\pi_R^\perp u_1\|_{L^2(\{x\in \Rm^6: |x|>R\})}^2.
\end{align*}
Here $\Pi_R^\perp$ is the orthogonal projection from $\dot{H}^1(\{x\in \Rm^4: |x|>R\})$ onto the complement space of $\hbox{span} \{|x|^{-2}\}$. While $\pi_R^\perp$ is the orthogonal projection from $L^2 (\{x\in \Rm^6: |x|>R\})$ onto the complement space of $\hbox{span} \{|x|^{-4}\}$. 

\subsection{Main idea}

According to \eqref{relationship radiation and channel} we may obtain exterior energy estimates conveniently from the radiation profiles $G_\pm$. Please note that $G_-$ and $G_+$ are not independent to each other. In fact the map $\mathbf{T}_+\circ \mathbf{T}_-^{-1} : G_- \rightarrow G_+$ is a bijective isometry. If we could find explicit expressions of the maps 
\begin{align*}
 &\mathbf{T}_+\circ \mathbf{T}_-^{-1} : G_- \rightarrow G_+;& &\mathbf{T}_-^{-1} : G_- \rightarrow (u_0,u_1);& &\mathbf{S}_L \circ \mathbf{T}_-^{-1}: G_- \rightarrow u;&
\end{align*}
then we would be able to 
\begin{itemize}
 \item[(a)] Understand how the asymptotic behaviour in one time direction determines the behaviour in the other time direction. This is known in the odd dimensional case, as shown (although not stated explicitly) in the proof of Proposition \ref{situation R0} by Duyckaerts-Kenig-Merle \cite{dkmnonradial}. In this work we will try to figure out the even dimensional case. 
 \item[(b)] Characterize (weakly) non-radiative solutions, especially in the radial case. We first determine all the radiation profiles $G_-$ so that 
 \[
  G_-(s,\theta)=G_+(s,\theta) = 0, \; s>R \quad \Leftrightarrow \lim_{t\rightarrow \pm \infty} \int_{|x|>|t|+R} |\nabla_{x,t} u(x,t)|^2 dx = 0;
 \] 
 then we may obtain all the non-radiative solutions (as well as their initial data) by applying the formula of $\mathbf{T}_-^{-1}$. In particular we prove that radial non-radiative solutions in the even dimension can be characterized in the same way as in the odd dimensions. 
 \item[(c)] Prove exterior energy estimates. We generalize the radial exterior energy estimates in $4$ and $6$ dimension to all even dimensions; we also prove a non-radial exterior energy estimate in the odd dimensions. In both applications (b) and (c) we follow the same roadmap: 
 \[
  \hbox{exterior energy} \leftrightarrow \hbox{radiation profile} \leftrightarrow \hbox{initial data}.
\]
\end{itemize}

\subsection{Main results} 

Now we give the statement of our results.  The details and proof can be found in subsequent sections. 

\begin{theorem} \label{main2}
Let $u$ be a finite-energy free wave with an even spatial dimension $d\geq 2$ and $G_+$, $G_-$ be the radiation profiles associated with $u$. Then we may give an explicit formula of the operator $\mathbf{T}_+ \circ \mathbf{T}_-^{-1}: G_- \rightarrow G_+$ 
\[
 G_+ (s, \theta) =(-1)^{d/2} \left(\mathcal{H} G_-\right) (-s,-\theta)
\]
Here $\mathcal{H}$ is the Hilbert transform in the first variable, i.e. 
\[
 \left(\mathcal{H} G_-\right) (-s,-\theta) = \hbox{p.v.}\, \frac{1}{\pi}\int_{-\infty}^\infty \frac{G_-(\tau, -\theta)}{(- s) - \tau} d\tau. 
\]
\end{theorem}

\begin{remark}
 A similar but simpler argument shows that if $d$ is odd, then $\mathbf{T}_+ \circ \mathbf{T}_-^{-1}: G_- \rightarrow G_+$ can be explicitly given by
 \[
   G_+ (s,\theta) = (-1)^{\frac{d-1}{2}} G_-(-s,-\theta). 
 \] 
This can also be verified by assuming that the initial data is smooth and compactly-supported, and considering the expression of $G_-$ and $G_+$ in terms of $(u_0,u_1)$ if $d$ is odd. Please refer to Duyckaerts-Kenig-Merle \cite{dkmnonradial}. Since we have $\mathcal{H}^2 = -1$. We may write the odd and even dimensions in a universal formula 
 \[
   G_+ (s,\theta) = ((-\mathcal{H})^{d-1} G_-) (-s,-\theta). 
 \] 
\end{remark}

\begin{remark}
 It has been proved in Section 3.2 of Duychaerts-Kenig-Merle \cite{oddtool} (in the language of Hankel and Laplace transforms) that the zero function is the only function $f \in L^2(\Rm)$ satisfying 
\begin{align*}
 &f(s) = 0, \; s>0;& &(\mathcal{H} f)(s) = 0,\; s<0.&
\end{align*}
It immediately follows that
\end{remark} 

\begin{corollary} 
Assume $d\geq 2$. Let $\Omega$ be a region in $\mathbb{S}^{d-1}$. If a finite-energy solution $u$ to homogenous linear wave equation satisfies 
\[
 \lim_{t\rightarrow \pm \infty} \int_{|t|}^\infty \int_{\pm \Omega} |\nabla_{t,x} u(r\theta, t)|^2 r^{d-1} d\theta dr = 0, 
\]
then we have
\[
 \lim_{t\rightarrow \pm \infty} \int_{0}^\infty \int_{\pm \Omega} |\nabla_{t,x} u(r\theta, t)|^2 r^{d-1} d\theta dr = 0.
\]
This is an angle-localized version of Corollary \ref{non radiative is zero}. 
\end{corollary}

\paragraph{Applications on channel of energy} By following the idea described above, we obtain the following results about the channel of energy. 

\begin{proposition} [Radial weakly non-radiative solutions] \label{main3}
 Let $d\geq 2$ be an integer and $R>0$ be a constant. If initial data $(u_0,u_1) \in \dot{H}^1 \times L^2$ are radial, then the corresponding solution to the homogeneous linear wave equation $u$ is $R$-weakly non-radiative, i.e. 
 \[
  \lim_{t\rightarrow \pm \infty} \int_{|x|>|t|+R} |\nabla_{t,x} u(x,t)|^2 dx = 0,
 \]
if and only if the restriction of $(u_0,u_1)$ in the region $\{x\in \Rm^d: |x|>R\}$ is contained in
\[
 P_{rad}(R) = \hbox{Span} \left\{(r^{2k_1-d},0), (0,r^{2k_2-d}): 1\leq k_1 \leq \left\lfloor \frac{d+1}{4} \right\rfloor, 1\leq k_2 \leq \left\lfloor \frac{d-1}{4} \right\rfloor\right\}
\]
Here the notation $\lfloor q\rfloor$ is the integer part of $q$. In particular, all radial $R$-weakly non-radiative solution in dimension $2$ are supported in $\{(x,t): |x|\leq |t|+R\}$. 
\end{proposition}
 
\begin{remark}
 If $d$ is odd, we have $\lfloor \frac{d+1}{4} \rfloor = \lfloor \frac{d+2}{4} \rfloor$ and $\lfloor \frac{d-1}{4} \rfloor = \lfloor \frac{d}{4} \rfloor$, thus our result here is the same as the already known result in odd dimension, as given in Theorem \ref{radialpi}. 
\end{remark} 
 
\begin{proposition} [Radial exterior estimates in even dimensions] \label{main4}
 Let $d =4k$ with $k \in \mathbb{N}$ and $R>0$. If initial data $u_0 \in \dot{H}^1(\Rm^d)$ are radial, then the corresponding solution $u$ to the homogenous linear wave equation with initial data $(u_0,0)$ satisfies 
\begin{align*}
 \lim_{t\rightarrow \infty} \int_{|x|>R+|t|} \left|\nabla u(x,t)\right|^2 dx = \lim_{t\rightarrow \infty} \int_{|x|>R+|t|} \left|u_t(x,t)\right|^2 dx \geq \frac{1}{4} \|\Pi_{Q_k(R)}^\perp u_0\|_{\dot{H}^1(\{x:|x|>R\})}^2.
\end{align*}
Here $\Pi_{Q_k(R)}^\perp$ is the orthogonal projection from $\dot{H}^1(\{x\in \Rm^d: |x|>R\})$ onto the complement of the $k$-dimensional linear space 
\[
 Q_k(R) = \hbox{Span} \left\{\frac{1}{|x|^{4k-2k_1}}: 1\leq k_1 \leq k\right\}.
\]
Similarly if the dimension $d=4k+2 \geq 2$ with $k\in \{0\}\cup \mathbb{N}$ and initial data $u_1\in L^2(\Rm^d)$ are radial, then  the corresponding solution $u$ to the homogenous linear wave equation with initial data $(0,u_1)$ satisfies 
\begin{align*}
 \lim_{t\rightarrow \infty} \int_{|x|>R+|t|} \left|\nabla u(x,t)\right|^2 dx = \lim_{t\rightarrow \infty} \int_{|x|>R+|t|} \left|u_t(x,t)\right|^2 dx \geq \frac{1}{4} \|\Pi_{Q'_k(R)}^\perp u_1\|_{L^2(\{x:|x|>R\})}^2.
\end{align*}
Here $\Pi_{Q'_k(R)}^\perp$ is the orthogonal projection from $L^2(\{x\in \Rm^d: |x|>R\})$ onto the complement of the $k$-dimensional linear space
\[
 Q'_k(R) = \hbox{Span} \left\{\frac{1}{|x|^{4k+2-2k_1}}: 1\leq k_1 \leq k\right\}.
\]
\end{proposition}
\begin{remark}
 Given $u_0 \in \dot{H}^1(\Rm^{4k})$ or $u_1\in L^2(\Rm^{4k+2})$, the orthogonal projection of $u_0$ or $u_1$ onto the finite dimensional space $Q_k(R)$ or $Q'_k(R)$ gradually vanishes as $R\rightarrow 0^+$. Therefore if we make $R\rightarrow 0^+$ in Proposition \ref{main4}, we immediately obtain \eqref{exterior full 4k} and \eqref{exterior full 4k plus 2}. 
\end{remark}
\begin{proposition} [Non-radial exterior energy estimates] \label{main1} 
Let $d\geq 3$ be an odd integer and $R>0$ be a constant. Then the following inequality holds for all $(u_0,u_1) \in \dot{H}^1 \times L^2(\Rm^d)$. 
 \[
  \sum_{\pm} \lim_{t\rightarrow \pm \infty} \int_{|x|>R+|t|} \left|\nabla_{t,x} \mathbf{S}_L(t)(u_0,u_1)(x,t)\right|^2 dx = \left\|\Pi_{P(R)}^\perp (u_0,u_1)\right\|_{\dot{H}^1 \times L^2(\Rm^d)}^2. 
 \]
 Here $\Pi_{P(R)}^\perp$ is the orthogonal projection from $\dot{H}^1 \times L^2(\Rm^d)$ onto the complement of the closed linear space
 \[
 P(R) = \left\{(u_0,u_1)\in \dot{H}^1\times L^2(\Rm^d): \lim_{t\rightarrow \pm \infty} \int_{|x|>R+|t|} \left|\nabla_{t,x} \mathbf{S}_L(u_0,u_1)(x,t)\right|^2 dx = 0\right\}. 
\]
\end{proposition}

\paragraph{Structure of this work} This work is organized as follows. In section 2 we deduce an explicit formula of $\mathbf{T}_-^{-1}$ in all dimensions. Then in Section 3 we prove the explicit formula of $\mathbf{T}_+ \circ \mathbf{T}_-^{-1}$ given in Theorem \ref{main2}. The rest of the paper is devoted to the applications in channel of energy. We characterize radial weakly non-radiative solutions in Section 4, prove radial exterior energy estimate for all even dimensions in Section 5 and finally give a short proof of non-radial exterior energy estimate in odd-dimensional space in Section 6. The appendix is concerned with Hilbert transform of a family of special functions, since the Hilbert transform is involved in the even dimensions. 
 
\paragraph{Notations} In this work we use the notation $C(d)$ for a nonzero constant determined solely by the dimension $d$. It may represent different constants in different places. This avoid the trouble of keeping track of the constants when unnecessary.  
 
\section{From Radiation Profile to Solution} 

Now we assume that $G_-(r,\theta)$ is smooth and compactly supported and give an explicit formula of the operator $\mathbf{T}_-^{-1}$. We consider the odd dimensions first. 

\subsection{Odd dimensions}

\begin{lemma} \label{T inverse odd}
Assume that $d\geq 3$ is odd. Let $G_-$ be a smooth function with $\hbox{supp}\, G_- \subset [-R,R] \times \mathbb{S}^{d-1}$. Then $(u_0,u_1) = \mathbf{T}_-^{-1} G_-$ satisfies 
\begin{align}
  u_0(x) & = \frac{1}{(2\pi)^\frac{d-1}{2}} \int_{\mathbb{S}^{d-1}}  G_-^{(\mu-1)} \left(x_0 \cdot \omega, \omega\right) d\omega\\
 u_1(x) & = \frac{1}{(2\pi)^\frac{d-1}{2}} \int_{\mathbb{S}^{d-1}}  G_-^{(\mu)} \left(x_0 \cdot \omega, \omega\right)  d\omega
\end{align}
Here the notation $G_-^{(k)}$ represents the partial derivative
\[
   G_-^{(k)} (s,\theta)= \frac{\partial^k G_-(s,\theta)}{\partial s^k}. 
\]
\end{lemma}
\begin{remark}
 This formula in 3-dimensional case was previously known. Please refer to Friedlander \cite{inverseradiation}, for example. 
\end{remark}
\begin{proof}
 Let $(u_0,u_1)=\mathbf{T}_-^{-1} G_-$ and $u = \mathbf{S}_L(u_0,u_1)$. Given a large time $t>0$, we choose approximated data $(v_{0,t}, v_{1,t}) \approx (u(\cdot,-t), u_t(\cdot,-t))$ as below: 
 \begin{align}
  &v_{1,t} (r\theta) = r^{-\mu} G_-(r-t, \theta),& &r>0, \, \theta \in \mathbb{S}^{d-1};&\\
  &v_{0,t} (r\theta) = -\chi(r/t) \int_{r}^{+\infty} r'^{-\mu} G_-(r'-t,\theta) dr',& &r>0, \, \theta \in \mathbb{S}^{d-1}.&
 \end{align}
Here $\mu = (d-1)/2$ and $\chi: \Rm \rightarrow [0,1]$ is a smooth center cut-off function satisfying
\[
 \chi(s) = \left\{\begin{array}{ll} 1, & s>1/2;\\ 0, & s< 1/4.\end{array}\right.
\]
It is clear that the data $(v_{0,t}, v_{1,t})$ are smooth and compactly-supported in $\{x: |x|<R+t\}$. A straight-forward calculation shows that 
\begin{align*}
 \int_0^\infty \int_{\mathbb{S}^{d-1}} \left|r^{\mu} u_{1,t}(r\theta) - G_-(r-t, \theta)\right|^2 d\theta dr &=0;\\
 \int_0^\infty \int_{\mathbb{S}^{d-1}} \left|r^{\mu} \partial_r v_{0,t}(r\theta) - G_-(r-t, \theta)\right|^2 d\theta dr & \lesssim 1/t;\\
 \int_{\Rm^d} \left(|\nabla v_{0,t}(x)|^2 - |\partial_r v_{0,t}(x)|^2\right) dx & \lesssim 1/t. 
\end{align*}
Thus by radiation field we have 
\[ 
 \lim_{t\rightarrow +\infty} \left\|(v_{0,t}, v_{1,t})-(u(\cdot,-t), u_t(\cdot,-t))\right\|_{\dot{H}^1 \times L^2(\Rm^d)} = 0.  
\]
Since the linear propagation operator $\mathbf{S}_L(t)$ preserves the $\dot{H}^1 \times L^2$ norm, we have 
\begin{equation} \label{convergence u01}
 \lim_{t\rightarrow +\infty} \left\|\begin{pmatrix} u_0\\ u_1\end{pmatrix}  - \mathbf{S}_L(t) \begin{pmatrix} v_{0,t}\\ v_{1,t}\end{pmatrix} \right\|_{\dot{H}^1 \times L^2(\Rm^d)} = 0.
\end{equation}
Next we use the explicit expression of linear propagation operator (see, for instance, Evans \cite{pdeevans}) and write $v = \mathbf{S}_L(v_0,v_1)$ in terms of $(v_0,v_1)$ when the initial are sufficiently smooth. 
\begin{align*}
  v(x,t) & = c_d \cdot \frac{\partial}{\partial t} \left(\frac{1}{t}\frac{\partial}{\partial t}\right)^{\mu -1}\left(t^{d-2} \int_{\mathbb{S}^{d-1}} v_0(x+ t\omega) d\omega\right)\\
 & \qquad + c_d \cdot  \left(\frac{1}{t}\frac{\partial}{\partial t}\right)^{\mu -1}\left(t^{d-2} \int_{\mathbb{S}^{d-1}} v_1(x + t\omega) d\omega\right)\\
 & = c_d t^\mu \int_{\mathbb{S}^{d-1}} \left[((w\cdot \nabla)^{\mu} v_0)(x + t \omega) + ((w\cdot \nabla)^{\mu-1} v_1)(x + t \omega)\right] d\omega \\
 & \qquad + \sum_{0\leq k<\mu}  A_{d,k} t^k \int_{\mathbb{S}^{d-1}} ((w\cdot \nabla)^{k} v_0)(x + t \omega) d\omega \\
 & \qquad + \sum_{0\leq k<\mu-1} B_{d,k} t^{k+1} \int_{\mathbb{S}^{d-1}} ((w\cdot \nabla)^{k} v_1)(x + t \omega) d\omega.
\end{align*}
Here $c_d = \frac{1}{2(2\pi)^{\frac{d-1}{2}}}$, $A_{d,k}$, $B_{d,k}$ (and $A'_{d,k}$, $B'_{d,k}$ below) are all constants. We may differentiate and obtain 
\begin{align*}
 v_t (x,t) & = c_d t^\mu \int_{\mathbb{S}^{d-1}} \left[((w\cdot \nabla)^{\mu+1} v_0)(x + t \omega) + ((w\cdot \nabla)^{\mu} v_1)(x + t \omega)\right] d\omega\\
 & \qquad + \sum_{1\leq k\leq \mu}  A'_{d,k} t^{k-1} \int_{\mathbb{S}^{d-1}} ((w\cdot \nabla)^{k} v_0)(x + t \omega) d\omega\\
 & \qquad + \sum_{0\leq k\leq \mu-1} B'_{d,k} t^{k} \int_{\mathbb{S}^{d-1}} ((w\cdot \nabla)^{k} v_1)(x + t \omega) d\omega.
\end{align*}
Now we plug in $(v_0,v_1) = (v_{0,t}, v_{1,t})$ with large time $t$. We observe that 
\begin{align} \label{upper bound add} 
 |(\omega \cdot \nabla)^k v_{j,t}(x+tw)|  \lesssim t^{-\mu}, \qquad j=0,1;\; k\geq 0; 
\end{align}
and ($r = |x+t\omega|$, $\theta = \frac{x+t\omega}{|x + t\omega|}$, $k=\mu-1,\mu$)
\begin{align*}
 ((w\cdot \nabla)^{k+1} v_{0,t})(x + t \omega) & = \left(\omega\cdot \theta\right)^{k+1} r^{-\mu} G_-^{(k)} \left(r-t, \theta\right) + O(t^{-\mu-1});\\
 ((w\cdot \nabla)^{k} v_{1,t})(x +t\omega) & = \left(\omega\cdot \theta\right)^{k} r^{-\mu} G_-^{(k)} \left(r-t, \theta\right) + O(t^{-\mu-1}).
\end{align*}
Thus 
\[
 \begin{pmatrix} w_{0,t} \\ w_{1,t} \end{pmatrix} = \mathbf{S}_L(t)  \begin{pmatrix} v_{0,t}\\ v_{1,t} \end{pmatrix}
\]
satisfies 
\begin{align*}
 w_{0,t} & = c_d \int_{\mathbb{S}^{d-1}} \left(\omega\cdot \theta \right)^{\mu-1} \left(1 +  \omega\cdot \theta\right) G_-^{(\mu-1)} \left(r-t, \theta\right)  d\omega + O(1/t);\\
 w_{1,t} &= c_d \int_{\mathbb{S}^{d-1}} \left(\omega\cdot \theta\right)^\mu \left(1 +  \omega\cdot \theta \right) G_-^{(\mu)} \left(r-t, \theta\right)  d\omega + O(1/t).
\end{align*}
Please note that the implicit constants in \eqref{upper bound add}, $O(t^{-\mu-1})$ and $O(1/t)$ above may depend on $x$ but remain to be uniformly bounded if $x$ is contained in a compact subset of $\Rm^d$. Next we observe the facts 
\begin{align*}
 &\theta(\omega) = \omega + O(1/t);& &r(\omega)-t = x_0 \cdot \omega + O(1/t);&
\end{align*}
and further simplify the formula
\begin{align*}
w_{0,t} & = 2c_d \int_{\mathbb{S}^{d-1}}  G_-^{(\mu-1)} \left(x \cdot \omega, \omega\right)  d\omega + O(1/t);\\
w_{1,t} &= 2c_d \int_{\mathbb{S}^{d-1}}  G_-^{(\mu)} \left(x \cdot \omega, \omega\right)  d\omega + O(1/t). 
\end{align*}
Finally we make $t\rightarrow +\infty$, utilize \eqref{convergence u01} and obtain 
\begin{align*}
u_0 & = 2c_d \int_{\mathbb{S}^{d-1}} G_-^{(\mu-1)} \left(x \cdot \omega, \omega\right)  d\omega;\\
u_1 & = 2c_d \int_{\mathbb{S}^{d-1}}  G_-^{(\mu)} \left(x \cdot \omega, \omega\right)  d\omega.
\end{align*}
We plug in the value of $c_d$ and finish the proof. 
\end{proof}

\begin{remark}
An explicit formula of the free wave $u = \mathbf{S}_L \mathbf{T}_-^{-1} G_-$ can be given by
 \[
  u(x,t) = \frac{1}{(2\pi)^\frac{d-1}{2}} \int_{\mathbb{S}^{d-1}} G_-^{(\mu-1)} (x \cdot \omega + t, \omega) d\omega.
 \]
This can be verified by a straight-forward calculation. One may check
\begin{itemize}
 \item The function $u$ above is a smooth solution to the homogenous linear wave equation;
 \item The initial data of $u$ are exactly those given in Lemma \ref{T inverse odd}. 
\end{itemize}
We may differentiate and obtain
 \begin{align*}
  u_t(x,t) & = \frac{1}{(2\pi)^\frac{d-1}{2}} \int_{\mathbb{S}^{d-1}} G_-^{(\mu)} (x \cdot \omega + t, \omega) d\omega;\\
  \nabla u(x,t) & = \frac{1}{(2\pi)^\frac{d-1}{2}} \int_{\mathbb{S}^{d-1}} G_-^{(\mu)} (x \cdot \omega + t, \omega) \,\omega \, d\omega
 \end{align*}
\end{remark}

\subsection{Even dimensions}

The formula of $\mathbf{T}_-^{-1}$ in even dimensions are a little more complicated. 
\begin{lemma} \label{T inverse even}
 Assume that $d\geq 2$ is even and $G_- \in C_0^\infty (\Rm\times \mathbb{S}^{d-1})$. Then the operator $\mathbf{T}_-^{-1}$ is given explicitly by
 \begin{align*}
 u_0 (x) & = \frac{\sqrt{2}}{(2\pi)^{d/2}} \cdot  \int_0^{\infty} \int_{\mathbb{S}^{d-1}} \frac{G_-^{(d/2-1)} \left(x\cdot \omega -\rho, \omega\right)}{\sqrt{\rho}} d\omega d\rho;\\
 u_1 (x) & =\frac{\sqrt{2}}{(2\pi)^{d/2}} \cdot  \int_0^{\infty} \int_{\mathbb{S}^{d-1}} \frac{G_-^{(d/2)} \left(x\cdot \omega -\rho, \omega\right)}{\sqrt{\rho}} d\omega d\rho.
\end{align*}
\end{lemma}
\begin{proof}
Without loss of generality let us assume $\hbox{supp}\, G_- \subset [-R_1,R_1]\times \mathbb{S}^{d-1}$. It is sufficient to show that given any $R_2>0$, the formula above holds for almost everywhere $x \in B(0,R_2)$. Let us use the notations $(u_0,u_1) = \mathbf{T}_-^{-1} (G_-)$ and $u = \mathbf{S}_L(u_0,u_1)$. We consider the approximated data 
 \begin{align}
  &v_{1,t} (r\theta) = r^{-\mu} G_-(r-t, \theta),& \\
   &v_{0,t} (r\theta) = -\chi(r/t) \int_{r}^{+\infty} r'^{-\mu} G_-(r'-t,\theta) dr',& &r>0, \, \theta \in \mathbb{S}^{d-1}.&
 \end{align}
 and 
 \[
  \begin{pmatrix} w_{0,t} \\ w_{1,t} \end{pmatrix} = \mathbf{S}_L(t) \begin{pmatrix} v_{0,t}\\ v_{1,t}\end{pmatrix}. 
 \]
 Here $\chi$ is the center cut-off function as given in the previous subsection. A basic calculation shows 
 \[ 
 \lim_{t\rightarrow +\infty} \left\|(v_{0,t}, v_{1,t})-(u(\cdot,-t), u_t(\cdot,-t))\right\|_{\dot{H}^1 \times L^2(\Rm^d)} = 0.  
\]
Thus 
 \begin{equation} \label{convergence of w to u}
 \lim_{t\rightarrow +\infty} \left\|(w_{0,t}, w_{1,t})-(u_0, u_1)\right\|_{\dot{H}^1 \times L^2(\Rm^d)} = 0.  
\end{equation}
Let us first recall the explicit formula of $v = \mathbf{S}_L (v_0,v_1)$ in the even dimensional case:
\begin{align*}
  v(x,t) & = c_d \cdot \frac{\partial}{\partial t} \left(\frac{1}{t}\frac{\partial}{\partial t}\right)^{\frac{d-2}{2}}\left(t^{d-1} \int_{\mathbb{B}^d} \frac{v_0(x+ ty)}{\sqrt{1-|y|^2}}dy\right)\\
 & \qquad + c_d \cdot  \left(\frac{1}{t}\frac{\partial}{\partial t}\right)^{\frac{d-2}{2}}\left(t^{d-1} \int_{\mathbb{B}^{d}} \frac{v_1(x + ty)}{\sqrt{1-|y|^2}} dy\right)\\
 & = c_d \cdot t^{d/2} \int_{\mathbb{B}^{d}} \frac{((y\cdot \nabla)^{d/2} v_0)(x + t y) + ((y\cdot \nabla)^{d/2-1} v_1)(x + ty)}{\sqrt{1-|y|^2}} dy \\
 & \qquad + \sum_{0\leq k<d/2}  A_{d,k} t^k \int_{\mathbb{B}^d} \frac{(y\cdot \nabla)^{k} v_0(x + t y)}{\sqrt{1-|y|^2}} dy \\
 & \qquad + \sum_{0\leq k<d/2-1} B_{d,k} t^{k+1} \int_{\mathbb{B}^{d}} \frac{(y\cdot \nabla)^{k} v_1(x + t y)}{\sqrt{1-|y|^2}} dy.
\end{align*}
Here $\mathbb{B}_d$ is the unit ball in $\Rm^d$ and $c_d = (2\pi)^{-d/2}$ is a constant. The notations $A_{d,k}$, $B_{d,k}$ (and $A'_{d,k}$, $B'_{d,k}$ below) represent constants. We differentiate and obtain 
 \begin{align*}
  v_t (x,t) & = c_d \cdot t^{d/2} \int_{\mathbb{B}^{d}} \frac{((y\cdot \nabla)^{d/2+1} v_0)(x + t y) + ((y\cdot \nabla)^{d/2} v_1)(x + t y)}{\sqrt{1-|y|^2}} dy \\
 & \qquad + \sum_{1\leq k\leq d/2}  A'_{d,k} t^{k-1} \int_{\mathbb{B}^d} \frac{(y\cdot \nabla)^{k} v_0(x + t y)}{\sqrt{1-|y|^2}} dy \\
 & \qquad + \sum_{0\leq k\leq d/2-1} B'_{d,k} t^{k} \int_{\mathbb{B}^{d}} \frac{(y\cdot \nabla)^{k} v_1(x + t y)}{\sqrt{1-|y|^2}} dy.
 \end{align*}
We plug in $(v_0,v_1) = (v_{0,t}, v_{1,t})$ and observe
\begin{align*}
 &\left|(y\cdot \nabla)^{k} v_{0,t}\right| \leq t^{-\frac{d-1}{2}},& &\left|(y\cdot \nabla)^{k} v_{1,t}\right| \leq t^{-\frac{d-1}{2}}.&
\end{align*}
This gives the approximation 
\begin{align*}
 w_{0,t} (x) & = c_d \cdot t^{d/2} \int_{\mathbb{B}^{d}} \frac{((y\cdot \nabla)^{d/2} v_{0,t})(r\theta) + ((y\cdot \nabla)^{d/2-1} v_{1,t})(r\theta)}{\sqrt{1-|y|^2}} dy + O(t^{-1/2});\\
 w_{1,t} (x) & = c_d \cdot t^{d/2} \int_{\mathbb{B}^{d}} \frac{((y\cdot \nabla)^{d/2+1} v_{0,t})(r\theta) + ((y\cdot \nabla)^{d/2} v_{1,t})(r\theta)}{\sqrt{1-|y|^2}} dy + O(t^{-1/2}).
\end{align*}
Here $r = |x + t y|$, $\theta = \frac{x + t y}{|x + t y|}$. Furthermore, we observe ($k = d/2, d/2-1$)
\begin{align*}
  ((y\cdot \nabla)^{k+1} v_{0,t})(r\theta) & = \left(y\cdot \theta\right)^{k+1} r^{-\frac{d-1}{2}} G_-^{(k)} \left(r-t, \theta\right) + O(t^{-\frac{d+1}{2}});\\
 ((y\cdot \nabla)^{k} v_{1,t})(r\theta) & = \left(y\cdot \theta\right)^{k} r^{-\frac{d-1}{2}} G_-^{(k)} \left(r-t, \theta\right) + O(t^{-\frac{d+1}{2}});
\end{align*}
and write 
\begin{align*}
 w_{0,t} (x) & = c_d \cdot t^{d/2} \int_{\mathbb{B}^{d}} \frac{\left(y\cdot \theta\right)^{d/2-1} (y\cdot \theta+1) r^{-\frac{d-1}{2}} G_-^{(d/2-1)} \left(r-t, \theta\right)}{\sqrt{1-|y|^2}} dy + O(t^{-1/2});\\
 w_{1,t} (x) & = c_d \cdot t^{d/2} \int_{\mathbb{B}^{d}} \frac{\left(y\cdot \theta\right)^{d/2} (y\cdot \theta+1) r^{-\frac{d-1}{2}} G_-^{(d/2)} \left(r-t, \theta\right)}{\sqrt{1-|y|^2}} dy + O(t^{-1/2}).
\end{align*}
Next we observe that if $|y| < 1-\frac{R_1+R_2}{t}$, then we have $r\leq t|y| +|x| < t-R_1$ thus $G_-^{(k)} \left(r-t, \theta\right) = 0$. As a result, we may restrict the domain of integral to 
\[ 
 \mathbb{B}_t = \left\{y\in \mathbb{B}^{d}: |y| \geq 1 - \frac{R_1+R_2}{t}\right\}.
\]
Because in the region we have 
\begin{align*}
 &\theta = \frac{y}{|y|} + O(1/t);& &y \cdot \theta = 1 + O(1/t);& &r = t + O(1).&
\end{align*}
We can simplify the formula 
\begin{align*}
 w_{0,t} (x) & = 2c_d \cdot t^{1/2} \int_{\mathbb{B}_{t}} \frac{G_-^{(d/2-1)} \left(r-t, y/|y|\right)}{\sqrt{1-|y|^2}} dy + O(t^{-1/2});\\
 w_{1,t} (x) & = 2c_d \cdot t^{1/2} \int_{\mathbb{B}_{t}} \frac{G_-^{(d/2)} \left(r-t, y/|y|\right)}{\sqrt{1-|y|^2}} dy + O(t^{-1/2}).
\end{align*}
Next we utilize the change of variables 
\begin{align*}
 y = (1-\rho/t) \omega, \qquad (\rho, \omega) \in (0,R_1+R_2) \times \mathbb{S}^{d-1},
\end{align*}
and the approximations
\begin{align*}
 &r-t = x\cdot \omega -\rho + O(1/t);& &\sqrt{1-|y|^2} = (1+O(1/t))\sqrt{2\rho/t};& &dy = (1+O(1/t))t^{-1} d\rho d\omega;&
\end{align*}
 to obtain
 \begin{align*}
 w_{0,t} (x) & = \sqrt{2} c_d \cdot  \int_0^{R_1+R_2} \int_{\mathbb{S}^{d-1}} \frac{G_-^{(d/2-1)} \left(x\cdot \omega -\rho, \omega\right)}{\sqrt{\rho}} d\omega d\rho + O(t^{-1/2});\\
 w_{1,t} (x) & =\sqrt{2} c_d \cdot  \int_0^{R_1+R_2} \int_{\mathbb{S}^{d-1}} \frac{G_-^{(d/2)} \left(x\cdot \omega -\rho, \omega\right)}{\sqrt{\rho}} d\omega d\rho + O(t^{-1/2}).
\end{align*}
Finally we recall \eqref{convergence of w to u}, make $t\rightarrow +\infty$ and conclude 
\begin{align*}
 u_0 (x) & = \sqrt{2} c_d \cdot  \int_0^{R_1+R_2} \int_{\mathbb{S}^{d-1}} \frac{G_-^{(d/2-1)} \left(x\cdot \omega -\rho, \omega\right)}{\sqrt{\rho}} d\omega d\rho;\\
 u_1 (x) & =\sqrt{2} c_d \cdot  \int_0^{R_1+R_2} \int_{\mathbb{S}^{d-1}} \frac{G_-^{(d/2)} \left(x\cdot \omega -\rho, \omega\right)}{\sqrt{\rho}} d\omega d\rho.
\end{align*}
This finishes the proof. 
\end{proof}
\begin{remark}
 If $d\geq 4$, the convergence \eqref{convergence of w to u} implies that $(w_{0,t}, w_{1,t})$ converges to $(u_0,u_1)$ in $L^{\frac{2d}{d-2}} \times L^2$ by Sobolev embedding. We may combine this convergence with the local uniform convergence given above to verify the identities above. This argument breaks down in dimension $2$. We given another argument below in dimension $2$. Given any test function $\varphi \in C_0^\infty (\Rm^2)$, integration by parts gives an identity
\[
 \int w_{0,t}(x) \nabla \varphi(x) dx = - \int \nabla w_{0,t} (x) \varphi (x) dx. 
\]
We recall the local uniform convergence of $w_{0,t}$ given above and the $L^2$ convergence of $\nabla w_{0,t} \rightarrow \nabla u_0$, then obtain
\[
 \int \left(\sqrt{2} c_2 \cdot  \int_0^{\infty} \int_{\mathbb{S}^{1}} \frac{G_- \left(x\cdot \omega -\rho, \omega\right)}{\sqrt{\rho}} d\omega d\rho\right)  \nabla \varphi(x) dx = - \int \nabla u_0(x) \varphi (x) dx.
\]
This finishes the proof. Finally the author would like to mention that we have
\[
 \lim_{|x|\rightarrow +\infty} \sqrt{2} c_2 \cdot  \int_0^{\infty} \int_{\mathbb{S}^{1}} \frac{G_- \left(x\cdot \omega -\rho, \omega\right)}{\sqrt{\rho}} d\omega d\rho = 0.
\] 
\end{remark}
\begin{corollary} \label{data other t}
If $G_- \in C_0^\infty (\Rm\times \mathbb{S}^{d-1})$, then $u = \mathbf{S}_L \mathbf{T}_-^{-1} (G_-)$ is given by
\[
 u(x,t) = \frac{\sqrt{2}}{(2\pi)^{d/2}} \int_0^{\infty} \int_{\mathbb{S}^{d-1}} \frac{G_-^{(d/2-1)} \left(x\cdot \omega -\rho+t, \omega\right)}{\sqrt{\rho}} d\omega d\rho.
\] 
Thus 
\[
 u_t (x,t) = \frac{\sqrt{2}}{(2\pi)^{d/2}} \int_0^{\infty} \int_{\mathbb{S}^{d-1}} \frac{G_-^{(d/2)} \left(x\cdot \omega -\rho+t, \omega\right)}{\sqrt{\rho}} d\omega d\rho.
\] 
\end{corollary} 
\begin{proof}
 A basic calculation shows that $u(x,t)$ solves the free wave equation with initial data given in Lemma \ref{T inverse even}. 
\end{proof}

\subsection{Universal formula} \label{sec: universal operators}
Now let us give a universal formula of $\mathbf{T}_-^{-1}$ for all dimensions. We first define two convolution operators ($1/\sqrt{\pi x}$ is understood as zero if $x<0$)
\begin{align*}
 &\mathcal{Q} f = \frac{1}{\sqrt{\pi x}} \ast f,& &\mathcal{Q}' f = \frac{1}{\sqrt{-\pi x}} \ast f.& 
\end{align*}
Their Fourier symbols are $\frac{1-i(\xi/|\xi|)}{2\sqrt{\pi |\xi|}}$ and $\frac{1+i(\xi/|\xi|)}{2\sqrt{\pi |\xi|}}$, respectively. Let us also use the notation $\mathcal{D} = d/dx$ and recall that its Fourier symbol is $2\pi i \xi$. A simple calculation of symbols shows
\begin{align} \label{operator relationship}
 &\mathcal{Q}^2 \mathcal{D}= 1;& &\mathcal{Q}'^2 \mathcal{D} = -1;& &\mathcal{Q} \mathcal{Q}' \mathcal{D} = \mathcal{H}.&
\end{align}
As a result, we may understand $\mathcal{Q}$ as $\mathcal{D}^{-1/2}$ and rewrite $u = \mathbf{S}_L \mathbf{T}_-^{-1} G_-$ in the form of 
\begin{align}
 u(x,t) & = \frac{1}{(2\pi)^{(d-1)/2}} \int_{\mathbb{S}^{d-1}} \left(\mathcal{Q} G_-^{(d/2-1)}\right)\left(x\cdot \omega +t, \omega\right) d\omega \nonumber \\
 & = \frac{1}{(2\pi)^{\mu}} \int_{\mathbb{S}^{d-1}} \mathcal{D}^{\mu-1} G_- \left(x\cdot \omega +t, \omega\right) d\omega.
\end{align}
Here $\mu = \frac{d-1}{2}$. This formula holds for both odd and even dimensions.

\section{Between Radiation Profiles}

In this section we give an explicit expression of the operator $\mathbf{T}_+ \circ \mathbf{T}_-^{-1}$ in the even dimension case, without the radial assumption. 

\begin{theorem}
Assume that $d\geq 2$ is an even integer. The operator $\mathbf{T}_+ \circ \mathbf{T}_-^{-1}$ can be explicitly given by the formula 
\[
 G_+(s,\theta) = \left(\mathbf{T}_+ \mathbf{T}_-^{-1} G_-\right) (s, \theta) =(-1)^{d/2} \left(\mathcal{H} G_-\right) (-s,-\theta)
\]
Here $\mathcal{H}$ is the Hilbert transform in the first variable, i.e. 
\[
 \left(\mathcal{H} G_-\right) (-s,-\theta) = \hbox{p.v.}\, \frac{1}{\pi}\int_{-\infty}^\infty \frac{G_-(\tau, -\theta)}{-\tau - s} d\tau. 
\]
\end{theorem}
\begin{proof}
 Since $\mathbf{T}_+ \circ \mathbf{T}_-^{-1}$ is a bijective isometry from $L^2(\Rm\times \mathbb{S}^{d-1})$ to itself. We only need to prove this formula for smooth and compactly supported data $G_-$. Without loss of generality let us assume $\hbox{supp}\, G_- \subset [-R_1,R_1]\times \mathbb{S}^{d-1}$. Let us also fix a positive constant $R_2>0$. If $(s,\theta) \in (-R_2,R_2) \times \mathbb{S}^{d-1}$, then we may apply Corollary \ref{data other t} and obtain 
\begin{align*} 
(t+s)^{\frac{d-1}{2}} \partial_t u((t+s)\theta,t) = \sqrt{2} c_d (t+s)^{\frac{d-1}{2}} \int_0^{\infty} \int_{\mathbb{S}^{d-1}} \frac{G_-^{(d/2)} \left((t+s)\theta \cdot \omega -\rho+t, \omega\right)}{\sqrt{\rho}} d\omega d\rho
\end{align*}
Let $M \gg R_1 +R_2 + 1$ be a large constant, we may split the integral above into two parts
\begin{align*}
 J_1 & =  \sqrt{2} c_d (t+s)^{\frac{d-1}{2}} \int_0^{\infty} \int_{\theta \cdot \omega < -1+M/t} \frac{G_-^{(d/2)} \left((t+s)\theta \cdot \omega -\rho+t, \omega\right)}{\sqrt{\rho}} d\omega d\rho;\\
 J_2 & =  \sqrt{2} c_d (t+s)^{\frac{d-1}{2}} \int_0^{\infty} \int_{\theta \cdot \omega \geq -1+M/t} \frac{G_-^{(d/2)} \left((t+s)\theta \cdot \omega -\rho+t, \omega\right)}{\sqrt{\rho}} d\omega d\rho;
\end{align*}
We may find an upper bound of $J_2$. In this region we have 
\[
 (t+s)\theta \cdot \omega +t \geq M - R_2  \quad \Rightarrow \quad G_-((t+s)\theta \cdot \omega -\rho +t) = 0, \hbox{if}\, \rho < M/2.
\]
Thus we may integrate by parts and obtain 
\[
 J_2 = C(d) (t+s)^{\frac{d-1}{2}} \int_0^{\infty} \int_{\theta \cdot \omega \geq -1+M/t} \frac{G_- \left((t+s)\theta \cdot \omega -\rho+t, \omega\right)}{\rho^{\frac{d+1}{2}}} d\omega d\rho
\]
Thus when $t$ is sufficiently large 
\begin{align*} 
 |J_2| & \lesssim t^{\frac{d-1}{2}} \int_{\theta \cdot \omega \geq -1+M/t} \int_{(t+s)\theta\cdot \omega +t -R_1}^{(t+s)\theta\cdot \omega +t +R_1} \frac{\left|G_- \left((t+s)\theta \cdot \omega -\rho+t, \omega\right)\right|}{\rho^{\frac{d+1}{2}}} d\rho d\omega\\
 & \lesssim  t^{\frac{d-1}{2}}  \int_{\theta \cdot \omega \geq -1+M/t} \int_{(t+s)\theta\cdot \omega +t -R_1}^{(t+s)\theta\cdot \omega +t +R_1} \frac{\left|G_- \left((t+s)\theta \cdot \omega -\rho+t, \omega\right)\right|}{\left|(t+s)\theta \cdot \omega +t\right|^{\frac{d+1}{2}}} d\rho d\omega \\
 & \lesssim t^{\frac{d-1}{2}} \int_{\theta \cdot \omega \geq -1+M/t} \frac{1}{|t\theta\cdot \omega +t|^{\frac{d+1}{2}}} d\omega\\
 & \lesssim 1/M.
\end{align*}
In the integral region of $J_1$, we have the approximation $\omega = -\theta + O(t^{-1/2})$. Thus we have 
\[
 J_1 = \sqrt{2} c_d t^{\frac{d-1}{2}} \int_0^{\infty} \int_{\theta \cdot \omega < -1+M/t} \frac{G_-^{(d/2)} \left((t+s)\theta \cdot \omega -\rho+t, -\theta\right)}{\sqrt{\rho}} d\omega d\rho + O(t^{-1/2}). 
\]
Next we utilize the change of variables (please refer to figure \ref{change of variables} for a geometrical meaning)
\begin{align*}
 \omega & = (-1+\rho'/t)\theta + \sqrt{(\rho'/t)(2-\rho'/t)} \varphi, \quad \rho' \in [0,M], \, \varphi \in \mathbb{S}^{d-2} = \{\varphi\in \mathbb{S}^{d-1}: \varphi \perp \theta\}. \\
  d\omega & = [1+O(1/t)](2\rho'/t)^{\frac{d}{2}-1}d\mathbb{S}^{d-2}(\varphi) \cdot \frac{d\rho'}{\sqrt{2\rho' t}} = [1+O(1/t)] (2\rho')^\frac{d-3}{2} t^{-\frac{d-1}{2}} d\mathbb{S}^{d-2}(\varphi) d\rho'.
\end{align*}
and obtain
\[
 J_1 = \frac{1}{2 \pi^{d/2}} \int_0^\infty \int_0^M \int_{\mathbb{S}^{d-2}} G_-^{(d/2)} (\rho' -\rho -s, -\theta) \rho'^{\frac{d-3}{2}} \rho^{-1/2} d\varphi  d\rho' d\rho + O(t^{-1/2}).
\]
\begin{figure}[h]
 \centering
 \includegraphics[scale=1.0]{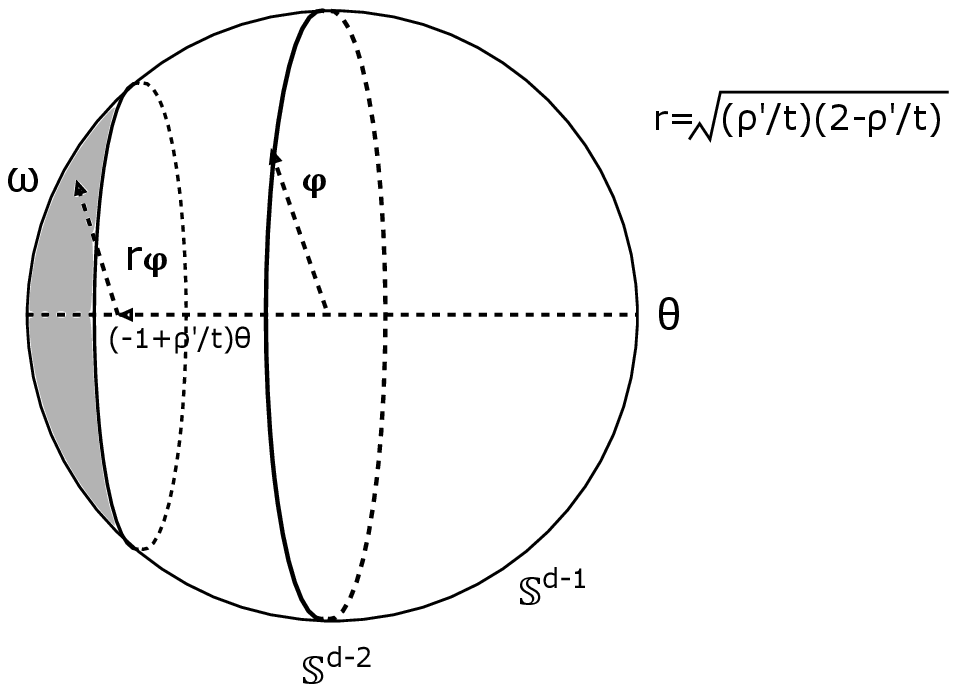}
 \caption{Change of variables} \label{change of variables}
\end{figure}
We observe that the integrand is independent of $\varphi$ and integrate by parts 
\[
 J_1 = \frac{(-1)^{d/2-1}}{\pi} \int_0^\infty \int_0^M \frac{G'_- (\rho' -\rho -s, -\theta)}{\sqrt{\rho \rho'}} d\rho' d\rho + O(t^{-1/2}). 
\]
We next change the variables $\tau = \rho' - \rho$, $\eta = \rho' + \rho$, and write
\begin{align*}
 J_1 & = \frac{(-1)^{d/2-1}}{\pi}  \int_{-\infty}^{M} \int_{|\tau|}^{2M-\tau} \frac{G'_-(\tau -s,-\theta)}{\sqrt{\eta^2 - \tau^2}} d\eta d\tau + O(t^{-1/2})\\
 & = \frac{(-1)^{d/2-1}}{\pi} \int_{-\infty}^M G'_-(\tau -s, -\theta) \left[\ln (2M-\tau+\sqrt{4M^2-4M\tau}) - \ln |\tau|\right] d\tau + O(t^{-1/2})\\
 & = \frac{(-1)^{d/2-1}}{\pi} \int_{-R_1-R_2}^{R_1+R_2} G'_-(\tau -s, -\theta) \left[\ln (2M-\tau+\sqrt{4M^2-4M\tau}) - \ln |\tau|\right] d\tau + O(t^{-1/2}).
\end{align*}
The integrals above can be split into two parts:
\begin{align*}
 I_1 & =  \int_{-R_1-R_2}^{R_1+R_2} G'_-(\tau -s, -\theta) \left[\ln (2M-\tau+\sqrt{4M^2-4M\tau}) \right] d\tau \\
 & = \int_{-R_1-R_2}^{R_1+R_2} G'_-(\tau -s, -\theta) \left[\ln (2M-\tau+\sqrt{4M^2-4M\tau}) - \ln (4M)\right] d\tau \\
 & = \int_{-R_1-R_2}^{R_1+R_2} G'_-(\tau -s, -\theta) O(1/M) d\tau = O(1/M);
\end{align*}
and 
\begin{align*}
 I_2 & = - \lim_{\eps\rightarrow 0^+} \int_{\eps<|\tau|<R_1+R_2} G'_-(\tau -s, -\theta) \ln |\tau| d\tau \\
 & =  \lim_{\eps\rightarrow 0^+} \int_{\eps<|\tau|<R_1+R_2} \frac{G_-(\tau -s, -\theta)}{\tau} d \tau \\
 & = -\pi (\mathcal{H} G_-) (-s,-\theta).
\end{align*}
In summary we have
\[
 J_1 = (-1)^{d/2} (\mathcal{H} G_-) (-s,-\theta) + O(1/M) + O(t^{-1/2}). 
\]
Now we may combine $J_1$ and $J_2$ 
\begin{align*}
 (t+s)^{\frac{d-1}{2}} \partial_t u((t+s)\theta,t) = (-1)^{d/2} (\mathcal{H} G_-) (-s, -\theta) + O(1/M) + O(t^{-1/2}).
\end{align*}
Because the implicit constants in $O$'s do not depend on $s\in [-R_2,R_2]$ or $\theta \in \mathbb{S}^{d-1}$, we may make $t\rightarrow +\infty$ then $M\rightarrow +\infty$ to conclude
\begin{align*}
 \lim_{t\rightarrow +\infty} \int_{-R_2}^{R_2} \int_{\mathbb{S}^{d-1}} \left|(t+s)^{\frac{d-1}{2}} \partial_t u((t+s)\theta,t) - (-1)^{d/2} (\mathcal{H} G_-) (-s, -\theta) \right|^2 d\theta ds  = 0.
\end{align*}
This finishes the proof.  
\end{proof}

\section{Radial Weakly Non-radiative Solutions}

In this section we prove Proposition \ref{main3}. First of all, we briefly show that any initial data in $P_{rad}(R)$ leads to a $R$-weakly non-radiative solution. By linearly we only need to consider the case $(u_0,u_1) = (r^{2k_1-d},0)$ or $(u_0,u_1) = (0, r^{2k_2-d})$. If $(u_0,u_1) = (r^{2k_1-d},0)$, then a basic calculation shows that if we choose $C_1, C_2, \cdots, C_{k_1-1}$ inductively, the solution
\[
 u_{k_1} (x,t) = \frac{1}{|x|^{d-2k_1}} + \frac{C_1 t^2}{|x|^{d-2k_1+2}} + \cdots + \frac{C_{k_1-1} t^{2k_1-2} }{|x|^{d-2}} 
\] 
solves the linear wave equation with initial data $(|x|^{2k_1-d},0)$ in the region $\Rm^d \setminus \{0\}$. By finite speed of propagation, we have
\[
 \mathbf{S}_L(u_0,u_1) (x,t) = u_{k_1} (x,t),\qquad |x|>R+|t|.
\] 
A simple calculation shows that this is indeed a non-radiative solution. The case $(u_0,u_1) = (0, r^{2k_2-d})$ can be dealt with in the same manner by considering the solution
\[
 u_{k_2}(x,t) =  \frac{t}{|x|^{d-2k_1}} + \frac{C_1 t^3}{|x|^{d-2k_1+2}} + \cdots + \frac{C_{k_2-1} t^{2k_1-1} }{|x|^{d-2}}.
\] 
Thus it is sufficient to show initial data of any non-radiative solution are contained in the space $P_{rad}(R)$. We first consider the odd dimensions. 

\subsection{Odd dimensions}

Assume that $u = \mathbf{S}_L (u_0,u_1)$ is a radial $R$-weakly non-radiative solution. Let $G_- = \mathbf{T}_-(u_0,u_1)$. By radial assumption $G_-$ is independent of the angle $\omega\in \mathbb{S}^{d-1}$. Let us first consider smooth functions $G_-$. We may calculate ($r>R$, $e_1 = (1,0,\cdots, 0) \in \Rm^d$)
\begin{align*}
 u_0(r e_1) =  (2\pi)^{-\mu} \int_{\mathbb{S}^{d-1}} G_-^{(\mu-1)} (r \omega_1) d\omega = \frac{\sigma_{d-2}}{(2\pi)^{\mu}}  \int_{-1}^{1} G_-^{(\mu-1)} (r \omega_1) (1-\omega_1^2)^{\mu-1} d\omega_1.
\end{align*}
Here $\omega_1$ is the first variable of $\Rm^d \supset \mathbb{S}^{d-1}$, $\sigma_{d-2}$ is the area of the sphere $\mathbb{S}^{d-2}$. We may integrate by parts and rescale 
\begin{align*}
 u_0 (r e_1) & = \frac{(-1)^{\mu-1} \sigma_{d-2}}{(2\pi)^{\mu} r^{\mu-1}} \int_{-1}^{1} G_-(r\omega_1) \left[\partial_{\omega_1}^{\mu-1} (1-\omega_1^2)^{\mu-1}\right] d\omega_1\\
 & = \sum_{k=0}^{\lfloor(\mu-1)/2\rfloor} \frac{A_{d,k}}{r^{\mu-1}}\int_{-1}^1 G_-(r\omega_1) \omega_1^{\mu-1-2k} d \omega _1\\
 & =  \sum_{k=0}^{\lfloor(d-3)/4\rfloor} \frac{A_{d,k}}{r^{d-2-2k}} \int_{-R}^R G_-(s) s^{\frac{d-3}{2}-2k} ds\\
 & =  \sum_{k=1}^{\lfloor(d+1)/4\rfloor} \frac{A_{d,k}}{r^{d-2k}} \int_{-R}^R G_-(s) s^{\frac{d+1}{2}-2k} ds
\end{align*}
Here $A_{d,k}$'s are nonzero constants. Similarly we have
\begin{align*}
 u_1(r e_1) & =  (2\pi)^{-\mu} \int_{\mathbb{S}^{d-1}} G_-^{(\mu)} (r \omega_1) d\omega \\
 & = \frac{\sigma_{d-2}}{(2\pi)^{\mu}} \int_{-1}^{1} G_-^{(\mu)} (r \omega_1) (1-\omega_1^2)^{\mu-1} d\omega_1\\
 & =  \frac{(-1)^\mu \sigma_{d-2}}{(2\pi)^{\mu} r^\mu} \int_{-1}^{1} G_-(r\omega_1) \left[\partial_{\omega_1}^{\mu} (1-\omega_1^2)^{\mu-1}\right] d\omega_1\\
 & = \sum_{k=0}^{\lfloor (\mu-2)/2 \rfloor} \frac{B_{d,k}}{r^\mu} \int_{-1}^1 G_-(r\omega_1) \omega_1^{\mu-2-2k} d \omega _1\\
 & = \sum_{k=1}^{\lfloor (d-1)/4 \rfloor} \frac{B_{d,k}}{r^{r-2k}} \int_{-R}^R G_-(s) s^{\frac{d-1}{2}-2k} ds. 
\end{align*}
Here $B_{d,k}$'s are nonzero constants. Since smooth functions are dense in $L^2([-R,R])$, we have
\begin{proposition} 
 There exist constants $\{A_{d,k}\}_{1\leq k \leq \lfloor(d+1)/4\rfloor}$, $\{B_{d,k}\}_{1\leq k \leq \lfloor(d-1)/4\rfloor}$, so that for any $G_- \in L^2(\Rm)$ supported in $[-R,R]$, the initial data $(u_0,u_1) = \mathbf{T}_-^{-1} G_-$ satisfy ($r>R$)
\begin{align*}
 u_0 (r) & = \sum_{k=1}^{\lfloor(d+1)/4\rfloor} \left(A_{d,k}  \int_{-R}^R G_-(s) s^{\frac{d+1}{2}-2k} ds\right) r^{-d+2k} ;\\
 u_1 (r) & = \sum_{k=1}^{\lfloor (d-1)/4 \rfloor} \left(B_{d,k}  \int_{-R}^R G_-(s) s^{\frac{d-1}{2}-2k} ds\right) r^{-d+2k}.
\end{align*}
\end{proposition}

\noindent This clearly shows that if $u = \mathbf{S}_L (u_0,u_1)$ is a radial $R$-weakly non-radiative solution, then $(u_0,u_1) \in P_{rad}(R)$. 

\subsection{Even dimensions}

The even dimensions involve Hilbert transform, thus are much more difficult to handle with. The general idea is the same. If the initial data $(u_0,u_1)$ are radial, then $G_\pm (s) = \mathbf{T}_\pm (u_0,u_1)$ is independent to the angle. We also have $G_+(s) = (-1)^{d/2} \mathcal{H} G_- (-s)$. Thus $\mathbf{S}_L (u_0,u_1)$ is $R$-weakly non-radiative if and only if $G_-$ is contained in the space
\[
 \mathcal{P}_{rad} = \{G_- \in L^2(\Rm): G_-(s)=0, s>R; \, (\mathcal{H}G_-)(s)=0, s<-R\}.
\]
Now recall the operators $\mathcal{Q}$, $\mathcal{Q}'$ and $\mathcal{D}$ defined in Subsection \ref{sec: universal operators}. We claim 
\begin{lemma} \label{structure of P radial}
 $\mathcal{Q}' \mathcal{P}_{rad} = H_0^{1/2} (-R,R) $. Here $\dot{H}_0^{1/2} (-R,R)$ is the completion of $C_0^\infty (-R,R)$ equipped with the $\dot{H}^{1/2}(\Rm)$ norm. 
\end{lemma}
\begin{proof}
 In order to avoid technical difficulties, we use an approximation technique. Given any $G_- \in \mathcal{P}_{rad}$, we may utilize a local smoothing kernel to generate a sequence $G_k$, so that 
\begin{itemize}
 \item[(a)] $G_k \in \mathcal{P}_{rad} (R+1/k)$; 
 \item[(b)] $G_k \in H^n (\Rm)$ for all $n\geq 0$ thus $G_k \in C^\infty(\Rm)$. 
 \item[(c)] $G_k$ converges to $G_-$ in $L^2(\Rm)$.  
\end{itemize} 
Let us consider the properties of the function $g_k = \mathcal{Q}' G_k \in C^\infty(\Rm)$. According to part (a), $G_k (s)= 0$ if $s>R+1/k$. We may use the convolution expression of $\mathcal{Q}'$ to obtain that $g_k$ vanishes in the interval $(R+1/k, +\infty)$. Similarly $g_k = \mathcal{Q} \mathcal{H} G_k$ vanishes in the interval $(-\infty, -R-1/k)$. We recall that $\mathcal{Q}' : L^2(\Rm) \rightarrow \dot{H}^{1/2}(\Rm)$ is an isometry up to a constant. Thus $g_k \rightarrow g = \mathcal{Q}' G_-$ in $\dot{H}^{1/2} (\Rm)$. This verifies $g \in \dot{H}_0^{1/2} (-R,R)$. We also need to show that given any $g \in \dot{H}_0^{1/2} (-R,R)$, then $\mathcal{Q}'^{-1} g \in \mathcal{P}_{rad}$. It is sufficient to consider $g \in C_0^\infty(-R,R)$ by smooth approximation. A simple calculation of Fourier symbols shows that $\mathcal{Q}'^{-1} = - \mathcal{Q}' \mathcal{D}$ and $\mathcal{H} \mathcal{Q}'^{-1} = \mathcal{Q} \mathcal{D}$. A combination of these identities with the convolution expressions of $\mathcal{Q}$ and $\mathcal{Q}'$ immediately verifies $\mathcal{Q}'^{-1} g \in \mathcal{P}_{rad}$. 
\end{proof}

\noindent We also need to use the following explicit formula of $\mathbf{T}_-$ for radial data 
\begin{lemma} \label{radial fast decay}
Assume $G \in C^\infty(\Rm)$ so that $|G(s)|\lesssim |s|^{-3/2}$ for $|s| \gg 1$. Then the corresponding radial free wave $u = \mathbf{S}_L \mathbf{T}_-^{-1} G$ satisfies
 \begin{equation} \label{radial fast decay formula}
  u(r,t) = C(d) \cdot r^{1-d/2} \int_{-1}^1 \mathcal{Q} G \left(r\omega_1 +t\right) P_d (w_1) (1-w_1^2)^{-1/2} d\omega_1.
 \end{equation} 
Here $P_d$ is an even or odd polynomial of degree $d/2-1$ defined by
\[
 \left(\frac{\partial}{\partial w_1}\right)^{\frac{d}{2}-1} (1-w_1^2)^\frac{d-3}{2} = P_d (w_1) (1-w_1^2)^{-1/2}. 
\] 
\end{lemma}
\begin{proof}  
If $G \in C_0^\infty (\Rm)$, we use the polar coordinates and integrate by parts:
\begin{align*}
 u(r,t) & = C(d) \int_0^{\infty} \int_{\mathbb{S}^{d-1}} \frac{G^{(d/2-1)} \left(r\omega_1 -\rho+t\right)}{\sqrt{\rho}} d\omega d\rho \\
 & =  C(d) \int_0^\infty \int_{-1}^1 \frac{G^{(d/2-1)} \left(r\omega_1 -\rho+t\right)}{\sqrt{\rho}} (1-w_1^2)^\frac{d-3}{2} d\omega_1 d\rho \\
 & = C(d) \cdot r^{1-d/2} \int_0^\infty \int_{-1}^1 \frac{G \left(r\omega_1 -\rho+t\right)}{\sqrt{\rho}} P_d (w_1) (1-w_1^2)^{-1/2} d\omega_1 d\rho\\
 & = C(d) \cdot r^{1-d/2} \int_{-1}^1 \mathcal{Q} G \left(r\omega_1 +t\right) P_d (w_1) (1-w_1^2)^{-1/2} d\omega_1.
\end{align*}
This verifies the formula if $G \in C_0^\infty (\Rm)$. In order to deal with profile $G$ without compact support, we use standard smooth cut-off techniques. More precisely, we may choose $G_k \in C_0^\infty (\Rm)$ so that $G_k \rightarrow G$ in $L^2(\Rm)$ and 
\[
 |G_k(s) - G(s)| \left\{\begin{array}{ll} = 0, & s<k;\\ \lesssim |s|^{-3/2} & s\geq k. \end{array}\right.
\]
Thus we have $\|\mathcal{Q} G - \mathcal{Q} G_k\|_{L^\infty} \lesssim 1/k$. This means we have the uniform convergence for all $(r,t)$ in any compact subset of $\Rm^+ \times \Rm$:
\begin{align*}
  u_k(r,t) & = \frac{C(d)}{r^{d/2-1}} \int_{-1}^1 \mathcal{Q} G_k \left(r\omega_1 +t\right) P_d(w_1) (1-w_1^2)^{-1/2} d\omega_1\\
   & \rightrightarrows \frac{C(d)}{r^{d/2-1}} \int_{-1}^1 \mathcal{Q} G \left(r\omega_1 +t\right) P_d(w_1) (1-w_1^2)^{-1/2} d\omega_1.
\end{align*}
Combining this with the convergence $u_k \rightarrow u$ in $\dot{H}^1$ we finish the proof. 
\end{proof}
\begin{remark}
 If $d\geq 4$ and $G\in L^2(\Rm)$, then formula \eqref{radial fast decay formula} still holds. This follows standard smooth approximation and/or cut-off techniques. Let $G_k \in C_0^\infty (\Rm)$ so that $G_k \rightarrow G$ in $L^2(\Rm)$. Thus $\mathcal{Q} G_k \rightarrow \mathcal {Q} G$ in $\dot{H}^{1/2}(\Rm)$. Finally we observe the fact $P_d(w_1) (1-w_1^2)^{-1/2} \in \dot{H}^{-1/2}(\Rm)$, obtain a locally uniform convergence $u_k(r,t)\rightarrow u(r,t)$ and conclude the proof. 
\end{remark}
\noindent Now we are ready to give an expression of $u = \mathbf{S}_L \mathbf{T}_-^{-1} G_-$ when $G_- \in \mathcal{P}_{rad}(R)$.
\begin{lemma} \label{expression even radial}
 Assume $G_- \in \mathcal{P}_{rad}$. Then the following identity holds
\[
 u(r,t) = \frac{C(d)}{r^{d/2}} \int_{-R}^{R} \mathcal{Q}' G_-(s) W_d\left(\frac{s-t}{r}\right) ds. 
\]
Here $W_d(\sigma)$ is the Hilbert transform (the functions below is understood as zero if $|w_1|>1$)
\[
 W_d(\sigma) \doteq \mathcal{H} \left(\frac{P_d(w_1)}{\sqrt{1-w_1^2}}\right) = \mathcal{H} \left[\left(\frac{d}{dw_1}\right)^{\frac{d}{2}-1}\left(1-w_1^2\right)^{\frac{d-3}{2}}\right];
\]
\end{lemma}
\begin{proof}
 By Lemma \ref{structure of P radial}, we have $\mathcal{Q}' G_- \in \dot{H}_0^{1/2}(-R,R)$. We claim that it is sufficient to consider the case $\mathcal{Q}' G_- \in C_0^\infty (-R,R)$. In fact, we may choose $G_k \in \mathcal{P}_{rad}(R)$ so that $\mathcal{Q}' G_k \in C_0^\infty (-R,R)$ so that 
\[
 \mathcal{Q}' G_k \rightarrow \mathcal{Q}' G_- \; \hbox{in}\; \dot{H}^{1/2}(-R,R) \quad \Leftrightarrow \quad G_k \rightarrow G_- \; \hbox{in}\; L^2(\Rm)
\]
Now we observe a few important facts: the embedding $\dot{H}_0^{1/2}(-R,R) \hookrightarrow L^p(-R,R)$ for all $1\leq p<+\infty$ and 
\[
 \frac{P_d(w_1)}{\sqrt{1-w_1^2}} \in L^p (\Rm) \quad \Rightarrow \quad W_d(\sigma) \in L^p(\Rm), \quad p\in (1,2).
\]
As a result, if the identity 
\[
 u_k(r,t) = \frac{C(d)}{r^{d/2}} \int_{-R}^{R} \mathcal{Q}' G_k(s) W_d\left(\frac{s-t}{r}\right) ds, \qquad k\geq 1
\]
holds, then we may make $k\rightarrow +\infty$ in the identity above and verify that a similar identity holds for $u$ and $G_-$. In fact the left hand side converges in the space $\dot{H}^1 (\Rm^d)$ for any given time $t$, while the right hand side converges uniformly for $(r,t)$ in any compact subset of $\Rm^+\times \Rm$. Now we assume $g = \mathcal{Q}' G_- \in C_0^\infty(-R,R)$. Then $G_- = \mathcal{Q}'^{-1} g = - \mathcal{Q}' \mathcal{D} g$ satisfies the assumption of Lemma \ref{radial fast decay}. As a result we have
\begin{align*}
  u(r,t) & = C(d) \cdot r^{1-d/2} \int_{-1}^1 \mathcal{Q} \mathcal{Q}' \mathcal{D} g \left(r\omega_1 +t\right) P_d (w_1) (1-w_1^2)^{-1/2} d\omega_1 \\
    & = C(d) \cdot r^{1-d/2} \int_{-1}^{1} \mathcal{H} g \left(r\omega_1 +t\right) P_d (w_1) (1-w_1^2)^{-1/2} d\omega_1\\
 & = \frac{C(d)}{r^{\frac{d}{2}-1}} \int_{-\infty}^\infty g (r\sigma +t) W_d (\sigma) d\sigma.
 \end{align*}
 Here we use the facts $\mathcal{Q} \mathcal{Q}' \mathcal{D} = \mathcal{H}$ and
\begin{align*}
&\int \mathcal{H} f \cdot \overline{\mathcal{H} g} dx = \int f\cdot \overline{g} dx& &\mathcal{H} (\mathcal{H} g(r\omega_1+t))(\sigma) = (\mathcal{H}^2 g) (r\sigma + t) = - g(r\sigma +t).&
\end{align*}
Finally we apply change of variables $s = r\sigma + t$, recall the support of $g$ and finish the proof. 
\end{proof}

\noindent Now let us consider the Hilbert transform $W_d$. The key observation is the following technical lemma. This result has probably been known for a long time, but we still give a brief proof in the Appendix for the purpose of completeness. 
\begin{lemma} \label{Hilbert special}
 Assume that $P(x)$ is a polynomial of degree $\kappa$. Let $W$ be the Hilbert transform 
 \[
  W = \mathcal{H}\left(\frac{P(x)}{\sqrt{1-x^2}}\right). 
 \]
Then $W(\sigma)$ is equal to a polynomial of degree $\kappa-1$ if $\sigma \in (-1,1)$. In particular, $W_2(\sigma)=0$ for $\sigma \in (-1,1)$; if $d\geq 4$, then the function $W_d(\sigma)$ is equal to an even or odd polynomial of degree $d/2-2$ in the interval $(-1,1)$. 
\end{lemma}

\paragraph{Proof of Proposition \ref{main4}} According to Lemma \ref{expression even radial}, we have already obtained 
\[
 u(r,t) = \frac{C(d)}{r^{d/2}} \int_{-R}^{R} \mathcal{Q}' G_-(s) W_d\left(\frac{s-t}{r}\right) ds.
\]
Here $\mathcal{Q}' G_- \in \dot{H}_0^{1/2}(-R,R) \hookrightarrow L^p (-R,R)$ for all $1<p<+\infty$. If we also have $r>|t|+R$, then 
\[
 \left|\frac{s-t}{r}\right| < 1, \qquad \forall s\in (-R,R).
\]
 If $d=2$, Lemma \ref{Hilbert special} immediately gives $u(r,t)\equiv 0$ if $r>|R|+t$ since we alway have $W_2(\frac{s-t}{r})=0$. In higher dimensional case $d\geq 4$, then Lemma \ref{Hilbert special} guarantees that 
\[
 W_d(s) = \sum_{l = 1}^{\lfloor d/4\rfloor} A_l s^{\frac{d}{2}-2l}, \qquad -1<s<1.
\]
is a polynomial. We plug this in the expression of $u$ and obtain
\begin{align} \label{expression u}
 u(r,t) = C(d) \sum_{l = 1}^{\lfloor d/4\rfloor} \frac{A_l}{r^{d-2l}} \int_{-R}^{R}  \mathcal{Q}' G_-(s) (s-t)^{\frac{d}{2}-2l} ds, \quad r>R +|t|.
\end{align}
This immediately gives $(u_0,u_1) \in P_{rad}(R)$.

\section{Exterior Energy Estimates of Even Dimensions}

In this section we prove Proposition \ref{main4}. It suffices to consider the case $d = 4k$. The proof of $d=4k+2$ are almost the same. Again we switch to the space of radiation profiles $G_- \in L^2 (\Rm \times \mathbb{S}^{d-1})$. We start by
\begin{lemma} \label{target of u0}
 The image of radial data in the form of $(u_0,0)$ can be characterized by
\begin{align*}
 \{\mathbf{T}_-(u_0,0): u_0 \in \dot{H}_{rad}^1(\Rm^d)\} & = \{G_-\in L^2(\Rm): \mathcal{H}G_- (-s) = - G_-(s)\} \\
 & = \left\{\frac{G(s)-\mathcal{H}G(-s)}{2}: G \in L^2(\Rm) \right\}
\end{align*} 
\end{lemma}
\begin{proof}
First of all, if $u_0 \in \dot{H}_{rad}^1 (\Rm^d)$, then free wave $u = \mathbf{S}_L (u_0,u_1)$ is radial and satisfies 
\begin{align*}
 &u(x,t) = u(x,-t);& &u_t(x,t) = - u_t(x,-t).&
\end{align*}
Therefore $G_-, G_+$ are radial, i.e. independent of $\omega$ and satisfy $G_+(s) = -G_-(s)$. We may apply Theorem \ref{main2} and obtain $G_+(s) = \mathcal{H} G_-(-s)$. As a result, $G_-$ satisfies the identity $\mathcal{H} G_- (-s) = - G_-(s)$. Next, let us assume $G_-$ satisfies this identity. Then we have 
\[
 G_-(s) = \frac{G_-(s)-\mathcal{H}G_-(-s)}{2} \in \left\{\frac{G(s)-\mathcal{H}G(-s)}{2}: G \in L^2(\Rm) \right\}. 
\]
Finally, if $G_-(s) = \frac{G(s) -\mathcal{H}G(-s)}{2}$, we show there exists $u_0 \in \dot{H}_{rad}^1 (\Rm^d)$, so that $G_- = \mathbf{T}_-(u_0,0)$. In fact, we consider radial initial data $(u_0,u_1) = \mathbf{T}_-^{-1} G$ and free wave $u = \mathbf{S}_L (u_0,u_1)$. We may reverse the time and obtain $u(x,-t) = \mathbf{S}_L (u_0,-u_1)(x,t)$. Thus 
\[
 \mathbf{T}_- (u_0,-u_1)(s) = -\mathbf{T}_+ (u_0,u_1)(s) = - \mathcal{H} G (-s)
\]
Therefore we have 
\[
 \mathbf{T}_-(2u_0,0) (s) = G(s) - \mathcal{H}G(-s) = 2 G_-(s)
\]
and complete the proof. 
\end{proof}
\noindent The key observation is the following
\begin{lemma} \label{Hilbert Laplace lemma}
 Given $g \in L^2(\Rm^+)$, there exists a function $G$ with $\|G\|_{L^2(\Rm)} \leq 2\|g\|_{L^2(\Rm^+)}$ so that 
\begin{align*}
 &G(s)-\mathcal{H}G(-s) = 2 g(s), \quad s>0. & &\left\|\frac{G(s)-\mathcal{H}G(-s)}{2}\right\|_{L^2(\Rm)} \leq \sqrt{2} \|g\|_{L^2(\Rm^+)}.&
\end{align*}
\end{lemma}
\begin{proof}
 Let us first find a function $G$ with $\|G\|_{L^2(\Rm)} \leq 2\|g\|_{L^2}$ so that
\[
 G(s) - \frac{G(s) + \mathcal{H} G(-s)}{2} = g(s), \quad s>0. 
\]
We define a linear bounded operator $\mathbf{T}$ from $L^2(\Rm^+)$ to itself\footnote{When we apply the Hilbert transform, we extend the domain of $G$ to $\Rm$ by assuming $G(s)=0$ if $s<0$.}
\[
 (\mathbf{T} G)(s) = \frac{G(s) + \mathcal{H} G(-s)}{2} = \frac{G(s)}{2} - \frac{1}{2\pi} \int_0^\infty \frac{G(\tau)}{s+\tau} d\tau, \quad s>0
\]
We may further rewrite it as 
\[
 \mathbf{T} G = \frac{G}{2} - \frac{1}{2\pi} \mathbf{L}^2 G. 
\]
Here $\mathbf{L}$ is the Laplace transform 
\[
 \mathbf{L} G(s) = \int_{0}^\infty G(\tau) e^{-s\tau} d\tau, 
\]
which is self-adjoint operator in $L^2(\Rm^+)$ with an operator norm $\sqrt{\pi}$. More details about the Laplace transform can be found in Lax \cite{lax}. As a result, we have 
\begin{align*}
 \|\mathbf{T} G\|_{L^2(\Rm^+)}^2 & = \frac{1}{4} \langle G- (1/\pi) \mathbf{L}^2 G, G- (1/\pi) \mathbf{L}^2 G \rangle\\
 & = \frac{1}{4} \|G\|_{L^2}^2 + \frac{1}{4\pi^2} \|\mathbf{L}^2 G\|_{L^2}^2 - \frac{1}{4\pi} \langle G, \mathbf{L}^2 G\rangle - \frac{1}{4\pi} \langle \mathbf{L}^2 G, G\rangle\\
 &  \leq \frac{1}{4} \|G\|_{L^2}^2 + \frac{1}{4\pi} \|\mathbf{L} G\|_{L^2}^2 - \frac{1}{2\pi} \langle \mathbf{L} G, \mathbf{L} G\rangle\\
 & =  \frac{1}{4} \|G\|_{L^2}^2 - \frac{1}{4\pi} \|\mathbf{L} G\|_{L^2}^2.
\end{align*}
Thus the operator norm of $\mathbf{T}$ is less or equal to $1/2$. This means that the function
\[
 G = \sum_{j=0}^\infty \mathbf{T}^j g \in L^2(\Rm^+)
\]
satisfies the equation $G - \mathbf{T} G = g$ and $\|G\|_{L^2(\Rm^+)} \leq 2\|g\|_{L^2(\Rm^+)}$. Finally we naturally extend the domain of $G$ to $\Rm$ by defining $G(s) = 0$ if $s<0$. We have 
\[
 \frac{G(s)-\mathcal{H}G(-s)}{2} = \left\{\begin{array}{ll} g(s), & s>0;\\ (-1/2) \mathcal{H} G(-s), & s<0. \end{array}\right.
\] 
Therefore we may find an upper bound of the $L^2$ norm
\begin{align*} 
 \left\|\frac{G(s)-\mathcal{H}G(-s)}{2}\right\|_{L^2(\Rm)}^2 & \leq \|g\|_{L^2(\Rm^+)}^2 + \frac{1}{4}\|\mathcal{H} G\|_{L^2(\Rm)}^2 \leq 2\|g\|_{L^2(\Rm^+)}^2.
\end{align*}
\end{proof}

\paragraph{Proof of Theorem \ref{main4}} Let $G_- = \mathbf{T}_- (u_0,0)$ and $g(s)$ be its cut-off version: 
\[
 g(s) = \left\{\begin{array}{ll} G_-(s), & s>R;\\ 0, & s<R. \end{array}\right. 
\]
Then radiation field implies that the free wave $u = \mathbf{S}_L (u_0,0)$ satisfies 
\begin{equation} \label{identity g negative}
\lim_{t\rightarrow -\infty} \int_{|x|>R+|t|} |\nabla u(x,t)|^2 dx =  \lim_{t\rightarrow -\infty} \int_{|x|>R+|t|} |u_t(x,t)|^2 dx = \sigma_{4k-1} \|g\|_{L^2(\Rm^+)}^2. 
\end{equation}
Here agian $\sigma_{4k-1}$ is the area of the sphere $\mathbb{S}^{4k-1}$. According to Lemma \ref{target of u0} and Lemma \ref{Hilbert Laplace lemma}, there exists a function $\tilde{u}_0 \in \dot{H}_{rad}^1 (\Rm^{4k})$, so that 
\begin{align*}
 &\mathbf{T}_- (\tilde{u}_0,0)(s) = g(s), \quad s>0; & &\|\tilde{u}_0\|_{\dot{H}^1(\Rm^{4k})}^2 \leq 4\sigma_{4k-1} \|g\|_{L^2(\Rm^+)}^2. &
\end{align*}
Therefore $\mathbf{T}_-(u_0-\tilde{u}_0, 0)$ vanishes if $s>R$. A combination of this fact with the time symmetry gives
\[
 \lim_{t\rightarrow \pm \infty} \int_{|x|>|t|+R} |\nabla_{t,x} \mathbf{S}_L (u_0-\tilde{u}_0, 0)(x,t)|^2 dx = 0.
\]
As a result, we may apply Proposition \ref{main3} and conclude $u_0 -\tilde{u}_0 \in Q_{k}(R)$. This means
\begin{align*}
 \left\|\Pi_{Q_k(R)}^\perp u_0\right\|_{\dot{H}^1(\{x: |x|>R\})}^2 \leq \|\tilde{u}_0\|_{\dot{H}^1(\{x: |x|>R\})}^2 & \leq 4\sigma_{4k-1} \|g\|_{L^2(\Rm^+)}^2 
\end{align*}
A combination of this inequality and identity \eqref{identity g negative} immediately verifies the conclusion of Proposition \ref{main4} in the negative time direction. The positive time direction follows the time symmetry. 
\section{Non-radial Exterior Energy Estimates}

In this section we give a short proof of Proposition \ref{main1}. We start by

\begin{lemma} 
 Let $d\geq 3$ be an odd integer. Then
 \begin{equation} \label{transformation PR}
 \sum_{\pm} \lim_{t\rightarrow \pm \infty} \int_{|x|>R+|t|} \left|\nabla_{t,x} \mathbf{S}_L (u_0,u_1)(x,t)\right|^2 dx = 2 \int_{|s|>R} \int_{\mathbb{S}^{d-1}} |\mathbf{T}_- (u_0,u_1)(s,\theta)|^2 d\theta ds. 
\end{equation}
In particular, we have (see \eqref{def of PR} for the definition of $P(R)$)
 \[
  \mathbf{T}_- (P(R)) = \mathcal{P}(R) \doteq \left\{G_-\in L^2(\Rm\times \mathbb{S}^{d-1}): \hbox{supp} \,G_- \subset [-R,R]\times \mathbb{S}^{d-1} \right\}. 
 \]
\end{lemma}
\begin{proof}
 Let $u$ be the solution of linear wave equation with initial data $(u_0,u_1)$. Then by radiation field (Theorem \ref{radiation}) we have 
 \begin{align*}
  \lim_{t\rightarrow -\infty} \int_{|x|>|t|+R} |\nabla_{t,x} u|^2 dx  & = 2 \int_{R}^\infty \int_{\mathbb{S}^{d-1}} |G_-(s,\theta)|^2 d\theta ds;\\
  \lim_{t\rightarrow -\infty} \int_{|x|<|t|-R} |\nabla_{t,x} u|^2 dx & = 2 \int_{-\infty}^{-R} \int_{\mathbb{S}^{d-1}} |G_-(s,\theta)|^2 d\theta ds.
 \end{align*}
 In addition, we may apply the energy conservation law, Proposition \ref{situation R0} and obtain
 \begin{align*}
  \lim_{t\rightarrow -\infty} \int_{|x|<|t|-R} |\nabla_{t,x} u|^2 dx & = \int_{\Rm^d} (|\nabla u_0|^2 +|u_1|^2) dx - \lim_{t\rightarrow -\infty} \int_{|x|>|t|-R} |\nabla_{t,x} u|^2 dx\\
  & = \lim_{t\rightarrow +\infty} \int_{|x|>t+R} |\nabla_{t,x} u|^2 dx. 
 \end{align*}
 Combining these identities we have 
\[
 \sum_{\pm} \lim_{t\rightarrow \pm \infty} \int_{|x|>R+|t|} \left|\nabla_{t,x} u(x,t)\right|^2 dx = 2 \int_{|s|>R} \int_{\mathbb{S}^{d-1}} |G_-(s,\theta)|^2 d\theta ds. 
\]
Finally $(u_0,u_1) \in P(R)$ is equivalent to saying
\[
 \int_{|s|>R} \int_{\mathbb{S}^{d-1}} |G_-(s,\theta)|^2 d\theta ds = 0,
\]
namely $\hbox{supp} \, G_- \subset [-R,R]\times \mathbb{S}^{d-1}$. This finishes the proof. 
\end{proof}

\noindent Now we are ready to prove Proposition \ref{main1}. Since $\sqrt{2} \mathbf{T}_-$ is a bijective isometry from $\dot{H}^1 \times L^2(\Rm^d)$ to $L^2(\Rm \times \mathbb{S}^{d-1})$. We have 
\[
 \mathbf{\Pi}_{P(R)}^\perp (u_0,u_1) = \mathbf{T}_-^{-1} \mathbf{\Pi}_{\mathbf{T}_- (P(R))}^\perp \mathbf{T}_- (u_0,u_1).
\]
We next use the expression of $\mathcal{P}(R) = \mathbf{T}_- (P(R))$: 
\begin{align*}
 \left\|\mathbf{\Pi}_{P(R)}^\perp (u_0,u_1)\right\|_{\dot{H}^1 \times L^2}^2 & = 2 \left\|\mathbf{\Pi}_{\mathcal{P}(R)}^\perp \mathbf{T}_- (u_0,u_1)\right\|_{L^2(\Rm\times \mathbb{S}^{d-1})}^2 \\
 &  = 2 \int_{|s|>R} \int_{\mathbb{S}^{d-1}} \left|\mathbf{T}(u_0,u_1)(s,\theta)\right|^2 d\theta ds. 
\end{align*}
Combining this with \eqref{transformation PR} we finish the proof. 

\section{Appendix}

 In this section we prove Lemma \ref{Hilbert special}. We first prove this lemma for two special cases, i.e. $P(x)=1$ and $P(x) = 1-x^2$. We start with $P(x)=1$. A straight forward calculate gives
\begin{align*}
 \pi W(s) & = \hbox{p.v.} \int_{-1}^1 \frac{(1-x^2)^{-1/2}}{s-x} dx \\
 & = \hbox{p.v.} \int_{-1}^1 \frac{(1-s^2)^{-1/2}}{s-x} dx + \int_{-1}^1 \frac{(1-x^2)^{-1/2}-(1-s^2)^{-1/2}}{s-x} dx \\
 & = (1-s^2)^{-1/2} \ln \left|\frac{1+s}{1-s}\right| + \int_{-1}^1 \frac{(1-s^2)-(1-x^2)}{(s-x)\sqrt{1-x^2} \sqrt{1-s^2} (\sqrt{1-x^2}+\sqrt{1-s^2})} dx\\
 & = (1-s^2)^{-1/2} \ln \left|\frac{1+s}{1-s}\right| + \frac{-s}{\sqrt{1-s^2}} \int_{-1}^1 \frac{1}{\sqrt{1-x^2}(\sqrt{1-x^2}+\sqrt{1-s^2})} dx. 
\end{align*}
Next we apply the change of variables $x = \frac{2z}{1+z^2}$. We have 
\begin{align*}
 &\sqrt{1-x^2} = \frac{1-z^2}{1+z^2}& &dx = \frac{2(1-z^2)}{(1+z^2)^2} dz&
\end{align*}
Thus 
\begin{align} 
 \int_{-1}^1 \frac{1}{\sqrt{1-x^2}(\sqrt{1-x^2}+\sqrt{1-s^2})} dx & = \int_{-1}^1 \frac{2dz}{1-z^2 + \sqrt{1-s^2} (1+z^2)} \nonumber\\
 & = \frac{2}{s} \int_{0}^1 \left(\frac{1}{\frac{1+\sqrt{1-s^2}}{s}-z} +\frac{1}{\frac{1+\sqrt{1-s^2}}{s}+z} \right) dz \nonumber \\
 & = \frac{2}{s} \ln \left|\frac{\frac{1+\sqrt{1-s^2}}{s}+1}{\frac{1+\sqrt{1-s^2}}{s}-1}\right| \nonumber \\
 & = \frac{1}{s} \ln \left|\frac{1+s}{1-s}\right| \label{special integral}
\end{align}
This immediately gives $W(x) = 0$. Next we consider the case $P(x) = 1-x^2$. In this case we calculate the Hilbert transform of $\sqrt{1-x^2}$
 \begin{align*}
 \pi W(s) & = \hbox{p.v.} \int_{-1}^1 \frac{\sqrt{1-x^2}}{s-x} dx \\
 & = \hbox{p.v.} \int_{-1}^1 \frac{\sqrt{1-s^2}}{s-x} dx + \int_{-1}^1 \frac{\sqrt{1-x^2}-\sqrt{1-s^2}}{s-x} dx \\
 & = \sqrt{1-s^2} \ln \left|\frac{1+s}{1-s}\right| + \int_{-1}^1 \frac{(1-x^2)-(1-s^2)}{(s-x)(\sqrt{1-x^2}+\sqrt{1-s^2})} dx\\
 & = \sqrt{1-s^2} \ln \left|\frac{1+s}{1-s}\right| + s\int_{-1}^1 \frac{1}{\sqrt{1-x^2}+\sqrt{1-s^2}} dx \\
 & = \sqrt{1-s^2} \ln \left|\frac{1+s}{1-s}\right| + \pi s + s\int_{-1}^1 \left(\frac{1}{\sqrt{1-x^2}+\sqrt{1-s^2}}-\frac{1}{\sqrt{1-x^2}}\right) dx\\
 & = \sqrt{1-s^2} \ln \left|\frac{1+s}{1-s}\right| + \pi s - s\sqrt{1-s^2} \int_{-1}^1 \frac{1}{\sqrt{1-x^2}(\sqrt{1-x^2}+\sqrt{1-s^2})} dx\\
 & = \pi s
\end{align*}
Here we use the integral \eqref{special integral} again. 
\paragraph{Induction} Now we are ready to prove Lemma \ref{Hilbert special} by an induction. It is clear that we only need to show the Hilbert transform of $f_\kappa(x) = x^\kappa (1-x^2)^{-1/2}$ is a polynomial of degree $\kappa-1$ in the interval $(-1,1)$. The cases of $\kappa =0, 2$ have been done. Now let us consider the case of $f_1(x) = x(1-x^2)^{-1/2}$. We observe that ($s\in (-1,1)$)
\[
 \mathcal{H} f_1 = \mathcal{H} \frac{d}{dx} \left(-\sqrt{1-x^2} \right) = - \frac{d}{ds} \mathcal{H} (\sqrt{1-x^2}) = -1.
\]
This prove the case $\kappa = 1$. Now let us assume that the cases $\kappa = 0,1,2,\cdots, n$ are done and consider the case $\kappa = n+1$. Here $n\geq 2$. We have 
\[
 x^{n+1} (1-x^2)^{-1/2} = - x^{n-1} (1-x^2)^{1/2} + x^{n-1} (1-x^2)^{-1/2}. 
\]
The Hilbert transform of the second term in the right hand side has been known to be a polynomial of degree $n-2$. Thus we only need to consider the first term. We have\footnote{Generally speaking, the derivative with respect to $s$ is in the weak sense. But since the derivative is known to be a polynomial in $(-1,1)$, we can integrate as usual.}
\begin{align*}
 \frac{d}{ds} \mathcal{H} \left(x^{n-1} (1-x^2)^{1/2}\right) & = \mathcal{H} \frac{d}{dx} \left(x^{n-1} (1-x^2)^{1/2}\right)\\
 & = \mathcal{H} \left\{\left[-n x^n + (n-1)x^{n-2}\right](1-x^2)^{-1/2} \right\}
\end{align*}
This is a polynomial of degree $n-1$ by induction hypothesis. A simple integration then finish the proof of case $\kappa = n+1$. 
\section*{Acknowledgement}
The second author is financially supported by National Natural Science Foundation of China Projects 12071339, 11771325.


\begin{thebibliography}{99}
 \bibitem{channeleven} R. C\^{o}te, C.E. Kenig and W. Schlag. {``Energy partition for linear radial wave equation.''} \textit{Mathematische Annalen} 358, 3-4(2014): 573-607.
 \bibitem{coratational} T. Duyckaerts, C.E. Kenig, Y. Martel and F. Merle. {Soliton resolution for critical co-rotational wave maps and radial cubic wave equation.} arXiv preprint 2103.01293.
 \bibitem{tkm1} T. Duyckaerts, C.E. Kenig, and F. Merle. {``Universality of blow-up profile for small radial type II blow-up solutions of the energy-critical wave equation.''} \textit{The Journal of the European Mathematical Society} 13, Issue 3(2011): 533-599.
  \bibitem{dkmnonradial} T. Duyckaerts, C. E. Kenig, and F. Merle. {``Universality of blow-up profile for small type II blow-up solutions of the energy-critical wave equation: the nonradial case''} \textit{The Journal of the European Mathematical Society} 14, Issue 5(2012): 1389-1454.
  \bibitem{se} T. Duyckaerts, C.E. Kenig, and F. Merle. {``Classification of radial solutions of the focusing, energy-critical wave equation.''} \textit{Cambridge Journal of Mathematics} 1(2013): 75-144.
 \bibitem{dkm2} T. Duyckaerts, C.E. Kenig, and F. Merle. {``Scattering for radial, bounded solutions of focusing supercritical wave equations.''} \textit{International Mathematics Research Notices} 2014:  224-258.
 \bibitem{dkm3} T. Duyckaerts, C.E. Kenig, and F. Merle. {``Scattering profile for global solutions of the energy-critical wave equation.''} \textit{Journal of European Mathematical Society} 21 (2019): 2117-2162.
 \bibitem{oddhigh} T. Duyckaerts, C. E. Kenig, and F. Merle. {``Soliton resolution for the critical wave equation with radial data in odd space dimensions.''}  arXiv preprint 1912.07664.
 \bibitem{oddtool} T. Duyckaerts, C. E. Kenig, and F. Merle. {``Decay estimates for nonradiative solutions of the energy-critical focusing wave equation.''} \textit{arXiv preprint} 1912.07665.
 \bibitem{pdeevans} L. C. Evans {``Partial Differential Equations, Second Edition.''} \textit{Graduate Studies in Mathematics} 19(2010), AMS, Providence.
\bibitem{radiation1} F. G. Friedlander. {``On the radiation field of pulse solutions of the wave equation.''}  \textit{Proceeding of the Royal Society Series A} 269 (1962): 53-65.
\bibitem{inverseradiation} F. G. Friedlander. {``An inverse problem for radiation fields.''} \textit{Proceeding of the London Mathematical Society} 27, no 3(1973): 551-576.
\bibitem{radiation2} F. G. Friedlander. {``Radiation fields and hyperbolic scattering theory.''} \textit{Mathematical Proceedings of Cambridge Philosophical  Society} 88(1980): 483-515.
\bibitem{fourierappli} G. B. Folland. {``Fourier analysis and its applications.''} \textit{The Wadsworth and Brooks/Cole mathematics series}, 1992, Pacific Grove, California. 
\bibitem{locad1} J. Ginibre, A. Soffer and G. Velo. {``The global Cauchy problem for the critical nonlinear wave equation''} \textit{Journal of Functional Analysis} 110(1992): 96-130.
 \bibitem{mg1} M. Grillakis. {``Regularity and asymptotic behaviour of the wave equation with critical nonlinearity.''} \textit{Annals of Mathematics} 132(1990): 485-509.
 \bibitem{mg2} M. Grillakis. {``Regularity for the wave equation with a critical nonlinearity.''} \textit{Communications on Pure and Applied Mathematics} 45(1992): 749-774.
\bibitem{loc1} L. Kapitanski. {``Weak and yet weaker solutions of semilinear wave equations''} \textit{Communications in Partial Differential Equations} 19(1994): 1629-1676.
\bibitem{katayamaradiation} S. Katayama. {``Asymptotic behavior for systems of nonlinear wave equations with multiple propagation speeds in three space dimensions.''} \textit{Journal of Differential Equations} 255(2013): 120-150.
 \bibitem{channel5d} C. E. Kenig, A. Lawrie, B. Liu and W. Schlag. {``Relaxation of wave maps exterior to a ball to harmonic maps for all data''} \textit{Geometric and Functional Analysis} 24(2014): 610-647.
 \bibitem{channel} C. E. Kenig, A. Lawrie, B. Liu and W. Schlag. {``Channels of energy for the linear radial wave equation.''} \textit{Advances in Mathematics}  285(2015): 877-936.
\bibitem{kenig} C. E. Kenig, and F. Merle. {``Global Well-posedness, scattering and blow-up for the energy critical focusing non-linear wave equation.''} \textit{Acta Mathematica} 201(2008): 147-212.
\bibitem{kenig1} C. E. Kenig, and F. Merle. {``Global well-posedness, scattering and blow-up for the energy critical, focusing, non-linear Schr\"{o}dinger equation in the radial case.''} \textit{Inventiones Mathematicae} 166(2006): 645-675.
 \bibitem{km} C. E. Kenig, and F. Merle. {``Nondispersive radial solutions to energy supercritical non-linear wave equations, with applications.''} \textit{American Journal of Mathematics} 133, No 4(2011): 1029-1065.
\bibitem{lax} P. D. Lax. {``Functional analysis.''} \textit{Pure and Applied Mathematics (New York)}, Wiley-Interscience, New York, 2002.
 \bibitem{ls} H. Lindblad, and C. Sogge. {``On existence and scattering with minimal regularity for semi-linear wave equations''} \textit{Journal of Functional Analysis} 130(1995): 357-426.
 \bibitem{enscatter1} K. Nakanishi. {``Unique global existence and asymptotic behaviour of solutions for wave equations with non-coercive critical nonlinearity.''} \textit{Communications in Partial Differential Equations} 24(1999): 185-221.
 \bibitem{enscatter2} K. Nakanishi. {``Scattering theory for nonlinear Klein-Gordon equations with Sobolev critical power.''} \textit{International Mathematics Research Notices} 1999, no.1: 31-60.
 \bibitem{ss1} J. Shatah, and M. Struwe. {``Regularity results for nonlinear wave equations''} \textit{Annals of Mathematics} 138(1993): 503-518.
  \bibitem{ss2} J. Shatah, and M. Struwe. {``Well-posedness in the energy space for semilinear wave equations with critical growth''} \textit{International Mathematics Research Notices} 7(1994): 303-309.
 \bibitem{shen2} R. Shen. {``On the energy subcritical, nonlinear wave equation in $\Rm^3$ with radial data''}  \textit{Analysis and PDE} 6(2013): 1929-1987.
\end{thebibliography}
\end{document}